# Intensive natural distribution as Bernoulli success ratio extension to continuous: enhanced Gaussian, continuous Poisson, and phenomena explanation

[1,*]Alessandro Felluga and [2]Stefano Tiziani

[1]ARPA Friuli Venezia Giulia, *Environm Protect Agcy* Lab of Udine, via Colugna 42, 33100, Italy

[2]Sanford-Burnham Medical Research Institute, 10901 North Torrey Pines Road, La Jolla, CA 92037, USA

[*]Corresponding author: Alessandro Felluga, alessandro.felluga@arpa.fvg.it

Stefano Tiziani, stiziani@sanfordburnham.org

**SUMMARY**

A new distribution called intensive natural distribution is introduced with the intent of merging statistics and empirical data. Based on the probability derived from the Bernoulli distribution, this method extended also Poisson distribution to continuous, preserving its skewness. Using this model, the Horwitz curve has been explained. The theoretical derivation of our method, which applies to every kind of measurements collected through sampling, is here supported by a mathematical demonstration and illustrated with several applications to real data collected from chemical and geotechnical fields. We compared the proposed intensive natural distribution to other widely-used frequency functions to test the robustness of the proposed method in fitting the histograms and the probability charts obtained from various intensive variables.

## 1. INTRODUCTION

Since the time of Gauss and Laplace, normal distribution has played an important role in both theoretical and statistical applications. Nowadays, the Gaussian model is still used as a first approximation for describing real situations in several fields; however, complex phenomena are hardly represented by the Gaussian distribution (Stingler (1977), Hill and Dixon (1982)). Statisticians have devoted intense efforts to the development of new parametric families in the attempt of depicting real data in a more realistic representation than using the normal model. Some of this literature includes flexible classes of skew distributions which based on the Gaussian model are capable of capturing skewness, tailweights, kurtosis and multimodality. Several families of skew-normal distributions have been introduced and subsequently developed by several authors including Pearson (1895), Johnson (1949), van Zwet (1964), Barndorff-Nielsen (1978), Azzalini (1985), Hoaglin (1986), Samorodnitsky and Taqqu (1994 ), Fernandez and Steel (1998), Azzalini and Capitanio (2003), Jones (2004), Azzalini (2005) and many others. Besides the skew-normal sampling distributions, alternative methods based on exponential power and skew-exponential power families have been proposed and discussed in detail (Azzalini (1986), Jones and Pewsey (2009)).

Although all these sampling distributions are flexible, widely applicable and capable of fitting the real data with a good approximation, they have a common limitation: the starting hypotheses are neglected therefore these constructed distributions fit the data without providing any interpretation about the phenomenon. “Why does the nature have a skewed probability distribution?” This question cannot be explained by using the abovementioned distributions, therefore a new statistic approach need to be explored starting from the beginning of the statistical model history.

In this paper, a novel statistic theory is introduced with the intent of merging statistics and empirical data using a comprehensive and complementary approach. Based on the probability derived from the Bernoulli distribution, this method is then extended to the continuous Poisson (1837) distribution similarly to the Gauss Law but avoiding critical approximations which compromise the interpretation of the phenomenon.

The logical path we adopted can be summarized as follows. The phenomenon needs to be explained therefore a model based on hypotheses is needed. The customary distribution which is based on the hypotheses for independent events is the binomial one, and this distribution can be skewed. Hence, the Poisson and the negative binomial distributions were derived. Starting from the probability theory, these distributions are discrete which intrinsically represent a too harsh approximation especially when continuous variables are observed. Gauss extended the Bernoulli distribution to continuous variables using the de Moivre-Laplace approximation which is based on the de Moivre-Stirling formula for large factorials. As a result, the Gaussian frequency function has maintained the strong Bernoulli hypotheses, the meaning and the application conditions. However, Gaussian curve is symmetric; the lack of skewness does not allow skewed data distribution interpretations. The Poisson and the negative binomial distributions are right side skewed but there is not a clear demonstration of their extension to continuous distributions. In addition, these curves are not derived from the initial hypotheses. The need of a phenomenological interpretation has then led statisticians to refer again to the binomial, Poisson or negative binomial distributions or their approximations and generalizations for data explanation (e.g. Ebneshahrashoob and Sobel (1990), Balasubramanian *et al.* (1993), Vellaisamy and Upadhye (2007), Chen *et al.* (2008), Malyshkina and Mannering (2008)).

All the above mentioned limitations suggest the necessity of a useful and simple continuous model based on solid hypotheses. In fact, from this “historical” research it is intuitive that intensive events have been implicitly described using homogeneous models. In reality, natural phenomena are dynamic and heterogeneous processes across time, space or both. The approximation of heterogeneity to homogeneity is the reason of the lack of skewness capability of the Gaussian frequency function. Pierre M. Gy (1992) studied extensively the heterogeneity producing an outstanding wide literature on the sampling theory. This author developed a theoretical model also with some parallelism with the symmetric model of the present paper but he focused mainly on the practical aspects of sampling on solid particulate matters without developing a general comprehensive theory, which is introduced in this paper.

## 2. DERIVATION OF THE MODEL FROM THE BASIS OF PROBABILITY

Let be M/N an intensive random variable defined as the ratio between two variables in which N be a very large number of events and M a large number of successes (favorable events) present in the N elements. For the central limit theorem, because N and M are both very large, then p, defined as M/N, is the single success probability. Assuming that the evaluation of all M successes in N events is not manageable, an estimation of p has to be performed by sampling. Note that each kind of estimation, including the instrumental analysis of a subsample, can be seen as a sampling/subsampling check. As a result, the variability of any type of measurements is due to the variability between the samples.

Since n is a sample collected within N events, inside n, m successes are observed. We express the bernoullian probability distribution in the sample of m success in n sampled elements, with p the single success probability to estimate as follows

$$pr(m,n) = \frac{n!}{m!(n-m)!} p^m q^{n-m}$$

where q = 1 – p is known to have expectation value E(m) = Np and variance var(m) = Npq. Assuming that n and m are large enough (for example 100), the factorial calculation is demanding also for modern hardwares. To minimize this computationally intensive factorial calculation, we apply the de Moivre-Stirling approximation for large factorials

$$pr(m,n) = \frac{\sqrt{2\pi n} \cdot n^n \cdot e^{-n}}{\sqrt{2\pi m} \cdot m^m \cdot e^{-m} \sqrt{2\pi(n-m)} \cdot (n-m)^{(n-m)} \cdot e^{-(n-m)}} p^m q^{n-m}$$

$$pr(m,n) = \sqrt{\frac{n}{2\pi m(n-m)}} \cdot \left(\frac{np}{m}\right)^m \cdot \left(\frac{nq}{m-n}\right)^{n-m}$$

where 0<m<n. In case m assumes values equal to 0 or n, the probability will not have any significance; while, increasing the number of n, the probability of m, for being equals to 0 or n, tends to zero. Assuming

$$\kappa = \left(\frac{np}{m}\right)^m \cdot \left(\frac{nq}{m-n}\right)^{n-m}$$

consequently

$$\begin{aligned} \mathit{ln}\,\kappa &= m\,\mathit{ln}\left(\frac{np}{m}\right) + (m-n)\mathit{ln}\left(\frac{nq}{m-n}\right) \\ &= m\,\mathit{ln}\left(1+\frac{\lambda}{m}\right) + (m-n)\mathit{ln}\left(1-\frac{\lambda}{m-n}\right) \end{aligned}$$

where λ = np-m; λ has the meaning of a distance from the expected value. But being n large and p finite, also m and n-m become large as np and nq, respectively. Because λ tends to zero, the use of natural logarithm expansion in McLaurin series is allowed

$$ln\,\kappa = m\left\{\left(\frac{\lambda}{m}\right) - \frac{1}{2}\left(\frac{\lambda}{m}\right)^2 + \frac{1}{3}\left(\frac{\lambda}{m}\right)^3 - \frac{1}{4}\left(\frac{\lambda}{m}\right)^4 + \ldots\right\}$$

$$+ (m-n)\left\{\left(\frac{\lambda}{n-m}\right) - \frac{1}{2}\left(\frac{\lambda}{n-m}\right)^2 + \frac{1}{3}\left(\frac{\lambda}{n-m}\right)^3 - \frac{1}{4}\left(\frac{\lambda}{n-m}\right)^4 + \ldots\right\}$$

$$= -\frac{\lambda^2}{2}\left(\frac{1}{m} + \frac{1}{n-m}\right) + \frac{\lambda^3}{3}\left(\frac{1}{m^2} + \frac{1}{(n-m)^2}\right) - \frac{\lambda^4}{4}\left(\frac{1}{m^3} + \frac{1}{(n-m)^3}\right) + \ldots$$

Because the ratio m/n tends to p which is a finite quantity, implicitly in case n diverges to infinite both m and (n – m) will diverge as well. In the last series, all terms tend to zero as fast as the power increases, therefore for large n, the sum is well approximated by the first term

$$ln\,\kappa \approx -\frac{\lambda^2}{2}\left(\frac{1}{m} + \frac{1}{n-m}\right) = -\frac{(Np-m)^2 n}{2m(n-m)}$$

then the approximated expression is obtained

$$pr(m,n) \approx \sqrt{\frac{n}{2\pi m(n-m)}} \cdot e^{-\frac{(Np-m)^2 n}{2m(n-m)}} \qquad (1)$$

So far the derivation is similar to the classic approximation of the binomial distribution of de Moivre-Laplace and similar to the Gy's Sampling Theory; the only adopted assumption is the large numbers of n and m. Starting from here, the discussion differs, because the customary models proceed using the following simplification

$$\frac{n}{m(n-m)} \approx Npq = \sigma^2 = const$$

for p not too close to its extreme values. Previously, the p to estimate was already extreme in case $p < 0.1\%$ or $p > 99.9\%$ when the differences to appreciate were expressed in "milli" units. Nowadays, owing to the increased sensitivity and reproducibility of the latest generation instruments, concentrations units are no longer only expressed in terms of "milli" units, but "micro", "pico" and "nano" units are becoming more realistically achievable and reproducible. As a result, the simplification used at the time of Gauss results too severe for modern high resolution analysis.

To maximize the empirical information obtained from modern instrumentation, we introduce the probability of the success ratio m/n. Considering

$$\mathrm{pr}\left(\frac{m}{n}\right)=\mathrm{pr}(m)\cdot\mathrm{pr}(n)=\mathrm{pr}(m)$$

we hypothesized that each sample has exactly n elemental unities, *ergo* pr(n) = 1.

Then

$$\mathrm{pr}\left(\frac{m}{n}\right)=\sqrt{\frac{n}{2\pi m(n-m)}}\cdot e^{-\frac{(Np-m)^2 n}{2m(n-m)}}=\frac{1}{\sqrt{2\pi n\frac{m}{n}\frac{n-m}{n}}}\cdot e^{-\frac{\left(p-\frac{m}{n}\right)^2}{2\frac{1}{n}\frac{m}{n}\frac{n-m}{n}}}$$

which is the probability of the ratio m/n; m/n is a discrete variable, which is changing by multiples of 1/n. Because n is quite large, the ratio m/n can be considered as a continuous positive variable $\chi$ in $\Re$ with $0 < \chi < 1$

$$\mathrm{pr}(\chi)=\frac{1}{\sqrt{2\pi n\chi(1-\chi)}}\cdot e^{-\frac{(p-\chi)^2}{2\frac{1}{n}\chi(1-\chi)}}$$

where $\chi$ is the success ratio of the sample. pr($\chi$) becomes infinitesimal as n diverges losing its significance. The frequency function remains finite because it is obtained dividing pr($\chi$) by the smallest possible interval 1/n which separates two contiguous values of $\chi$. As a result, for $0 < \chi < 1$ and n very large, the adimensional form can be rewritten as

$$f(\chi)=\frac{1}{\sqrt{2\pi\frac{1}{n}\chi(1-\chi)}}\cdot e^{-\frac{(p-\chi)^2}{2\frac{1}{n}\chi(1-\chi)}} \qquad (2)$$

which is the frequency function of the bernoullian success ratio $\chi$ (where $\chi$ is the success fraction in a sample of n elemental units of N and p is the global mean value of the ratio in the batch to estimate) extended to continuous $\chi$ values.

Assuming that the intensive quantity is measured with a new variable defined as the ratio between different unities for m and for n, the new intensive variable is more useful than the adimensional $\chi$ of equation (2). Introducing the constant $k \propto 1/n$ with n in the new units, k can assume all positive values ($0 < k < +\infty$) without modifying the equation (2). By substitution of the variable $\chi$, ranging from 0 to 1, with a variable x, between 0 and a certain maximum value u with $1 < u < +\infty$, also the value being estimated can be changed from p to $\mu$ (for analogy with the Gaussian model). Then equation (2) can be rearranged using $x = u\chi$, and consequently $\mu = up$. So being

$$f(x) = f(\chi)\left|\frac{d\chi}{dx}\right| = f(\chi)\left|\frac{1}{u}\right|$$

the frequency function of x, for x, $\mu \in \Re$ and $0 < x, \mu < u$, can be expressed as

$$f(x) = \frac{1}{\sqrt{2\pi k x(u - x)}} e^{-\frac{(\mu - x)^2}{2kx(u-x)}} \qquad (3)$$

This is the frequency function of the sample value x which derives from a batch with a $\mu$ mean value to estimate. Note that a third variable u is now present, which is the maximum possible value of x.

This curve can be symmetric, or lightly or heavily skewed. **Figure 1** shows three types of approximations which cover a comprehensive range of applications of this model to empirical data: the (Gaussian) normal-homogeneous approximation, the unlikely approximation and the likely approximations.

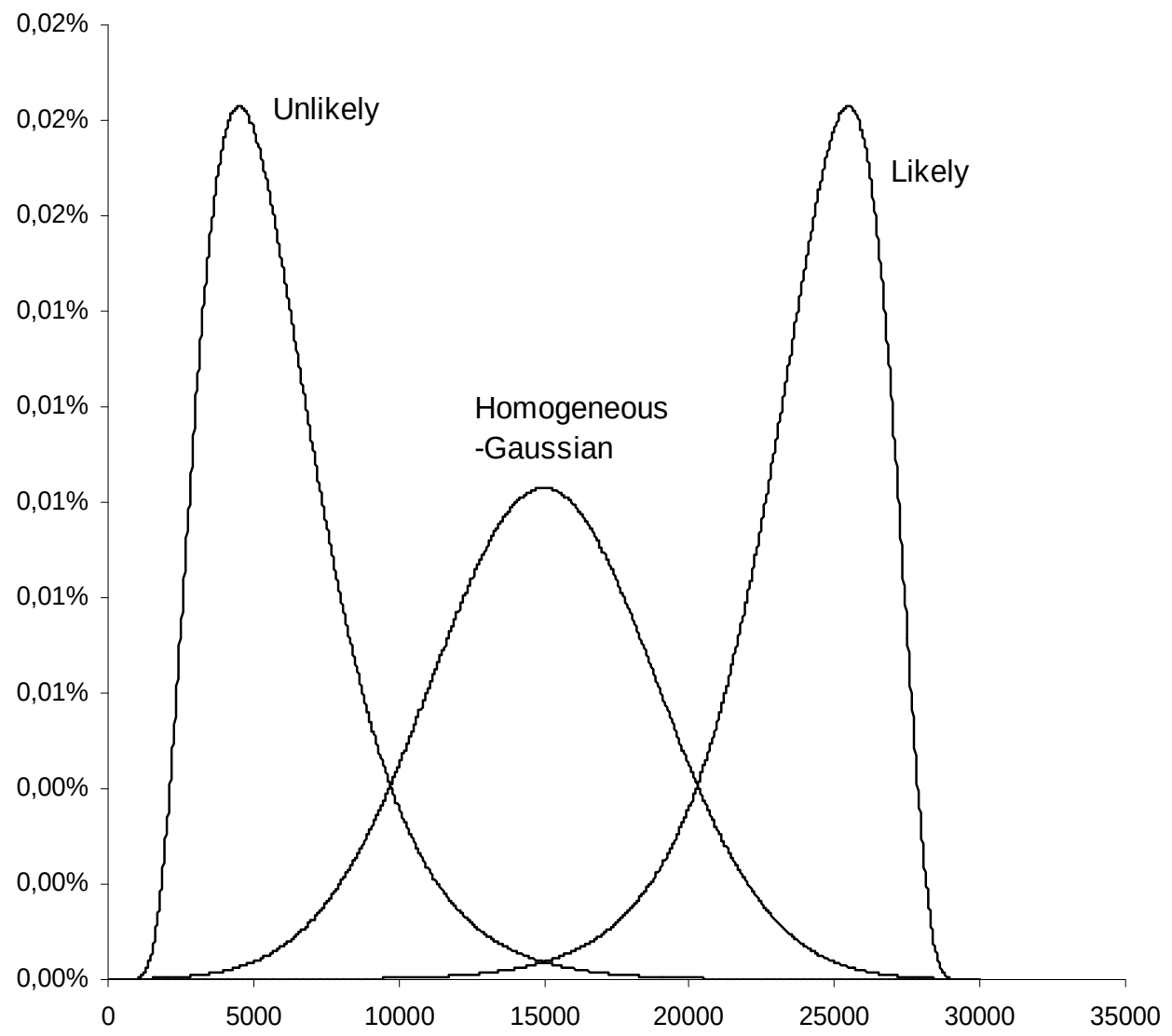

**Figure 1**. Three approximations of the intensive natural distribution equation (3) with the parameters with k = 1000 and m = 5000 for the Unlikely, m = 15000 for the Homogeneous-Gaussian and m = 25000 and u = 30000 for the Likely approximation.

## 3. APPROXIMATED FORMULAE OF THE GENERIC FORM

*3.1 The homogeneous intensive: the (enhanced) Gaussian curve.*

This assumption can be done when the batch to evaluate is homogeneous. For each sample collected, homogeneity means that the value of x is close to the batch mean value μ. So, under normal-homogeneous assumption, we can assume $x \approx \mu$, therefore equation (3) can be written as

$$f(x) = \frac{1}{\sqrt{2\pi k\mu(u-\mu)}} e^{-\frac{(\mu - x)^2}{2k\mu(u-\mu)}} \tag{4}$$

For $\sigma^2 = k\mu(u-\mu)$, equation (4) becomes exactly the Gaussian frequency function. The same result is obtained assuming *a priori* $\frac{n}{m(n-m)} \approx Npq = \sigma^2 = \text{const}$ as in the classical Gaussian model derivation showed before; therefore this assumption equals to a homogeneous batch approximations.

Equation (4) is an enhanced or "explicit" Gaussian distribution, a more comprehensive expression because it rationalizes heteroscedasticity (variance changing with the mean). It also prevents the occurrence of values, negative or over the maximum (u) because $\lim_{x \to 0^+} f(x) = \lim_{x \to u^-} f(x) = 0$. A plot of a function at low values is showed in **Figure 2**.

Equation (4) has an expectation value E(x) = μ and var(x) = kμ(u-μ) and it has all the well-known Gaussian frequency function properties.

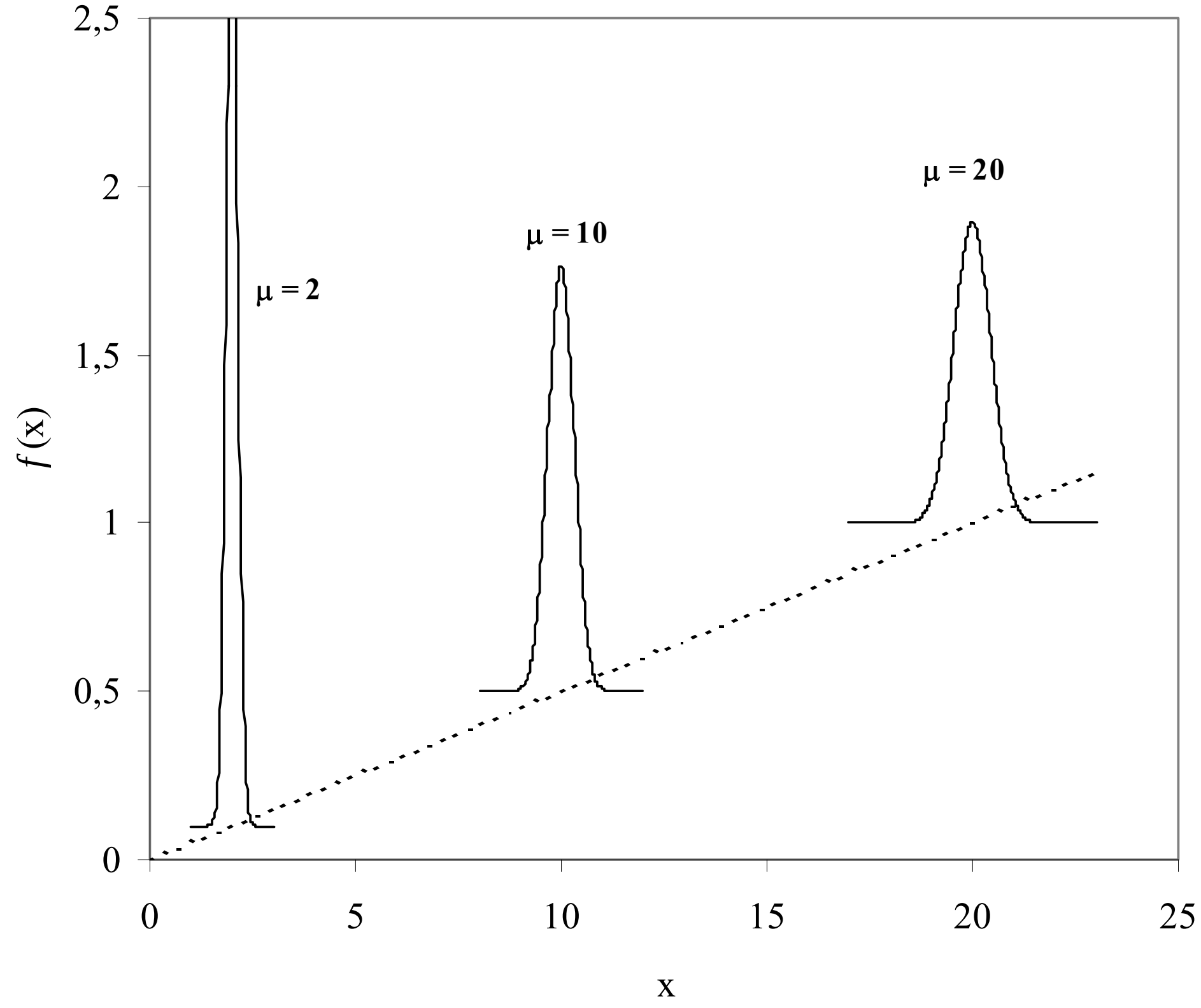

**Figure 2**. Enhanced Gaussian curve (homogeneous approximation, u >> μ, k = 0.01).

3.1.1 *Horwitz curve explanation.*

Horwitz *et al.* (1980) derived the so called Horwitz curve experimentally fitting a large number of coefficients of variation versus the general concentration levels of analytical measurements derived from collaborative trials studies. Horwitz found that the coefficient variation (CV) is a function of the concentration level expressed in adimensional mass/mass ratio p, regardless of the kind of analyte, matrix, instrumental technique, and other variables.

Considering that any single result in mass fraction $\chi$ is determined from a different sample and that all samples derive from a highly homogeneous material with a $\mu$ mean value, we can apply the intensive natural model for between-samples distribution with the homogeneous approximation. From the adimensional equation (2), considering homogeneity $\chi \approx p$ where p is the assigned value or target value of the collaborative exercise, this equation can be written as

$$f(\chi) = \frac{1}{\sqrt{2\pi \frac{1}{n} p(1-p)}} \cdot e^{-\frac{(p-\chi)^2}{2\frac{1}{n}p(1-p)}}$$

$\chi$ are the collected data and $E(\chi) = p$ and $Var(\chi) = p(1-p)/n$ are estimated to obtain a coefficient of variation, which is theoretically

$$CV = \frac{\sqrt{Var(\chi)}}{E(\chi)} = \frac{\sqrt{\frac{1}{n}p(1-p)}}{p} = \sqrt{\frac{1-p}{np}} \approx (np)^{-0.5}$$

When the value of p to estimate is low, *i.e.* ranging from $10^{-3}$ to $10^{-10}$ in mass ratio, p can be neglected ($p << 1$). This equation has a shape similar to the Horwitz curve. Hall and Selinger (1989) argued that the Horwitz curve can be perfectly reproduced starting from a CV formula $CV = np^{-0.5}$. The authors postulated that, starting from generic binomial trials and invoking the Zipf's law, n can be expressed as a function of p and for example an estimated $n(p) = 2500p^{0.3}$ (Hall and Selinger (1989)). The final equation is CV = $0.02p^{-0.15}$ which is similar to the empirical Horwitz equation. The present intensive natural model gives the necessary support to the Hall and Selinger's postulate completing the demonstration of the empirical Horwitz equation. Inferring that the Horwitz equation is obtained from both the intensive natural model distribution and the Zipf's principle of the least effort, the Horwitz curve represents the acceptable (minimum) performance expressed as CV between samples achieved with the minimum required effort. An extra effort at low level concentrations is achieved to reduce the naturally increasing of the variability at the lowest level concentration p; in fact the lower the p value the higher the between-samples variability.

*3.2 The unlikely intensive: a right skewed model.*

Consider now the unlikely event: the probability p is little ($p << 1$) and consequently the $\chi$ values having not zero probability are $\chi << 1$. This is like considering the generic form of equation (3) where $x << u$. In this approximation and re-labeling k, the product ku, the frequency function (3) is given by

$$f(x) = \frac{1}{\sqrt{2\pi k x}} \cdot e^{-\frac{(\mu - x)^2}{2kx}} \tag{5}$$

The normalization of the function (5) has been confirmed theoretically (PROOF 1). The batch mean value $\mu$ and the dispersion value k are (PROOF 2 and PROOF 3):

$$E\left(\frac{1}{x}\right) = \frac{1}{\mu} \qquad\qquad k = E(x) - \mu$$

Interestingly, the population mean $\mu$ is estimated by the harmonic mean of the samples values (and not by the simple mean) and the dispersion parameter k is the difference between the population mean $\mu$ and the samples mean E(x). Note that the mean value between samples is $E(x) = k + \mu$ and every time there is a systematic overestimation of the overall mean value μ due to the right tail. Increasing the sample size n, k decreases (as $k \propto 1/n$), therefore E(x) tends to μ as predicted by the central limit theorem.

Using the *maximum log-likehood* method the same expressions for the parameters μ and k would be

$$\frac{1}{\mu} = \frac{1}{N}\sum_{i=1}^{N}\frac{1}{x_i} \qquad \text{and} \qquad k = \frac{1}{N}\sum_{i=1}^{N}\frac{(x_i-\mu)^2}{x_i} = \frac{1}{N}\sum_{i=1}^{N}x_i - \mu$$

hence

$$\frac{1}{\mu} = E\left(\frac{1}{x}\right) \qquad \text{and} \qquad k = E(x) - \mu\,.$$

This systematic error in the overestimation of the mean E(x) adds to the random error of the mean of Student's theory (Student (1908)); however, the random error can be reduced by increasing the number of samples to mediate the systematic error. Therefore the overestimation is constant independently of the number of samples and it can be reduced only increasing the sample volume. The variance var(x) is (PROOF 4):

$$\operatorname{var}(x) = k(2k + \mu) \qquad \text{or} \qquad \operatorname{var}(x) = k^2 + kE(x)$$

Note that the variance is given by a fixed contribution $k^2$ which is due to the heterogeneity, therefore a sampling contribution and a variable contribution are related to both heterogeneity and quantity level.

The moment generating function (PROOF 5) for the unlikely intensive natural distribution is

$$M_x(t) = \frac{1}{\sqrt{1-2kt}}\, e^{\frac{\mu}{k}\left(1-\sqrt{1-2kt}\right)}$$

By sequential derivation of $M_x(t)$ in $t = 0$ all moments can be obtained; the calculated E(x) and Var(x) with this method confirmed the above written expressions.

The analogous characteristic function has been derived. It can be shown that if S is the sum of N random variables unlikely distributed, S is not unlikely distributed. This property, typical Gaussian and widely used for ambiguous composition, is here lost; the same was observed in skew-normal distributions (Dominguez-Molina and Rocha-Arteaga (2007), Kozubowski and Nolan (2008)).

Equation (5) is also achieved starting from the Poisson distribution with analogous steps passing to a continuous intensive variable; hence the unlikely intensive natural distribution can be expressed as a Poisson extended to continuous (PROOF 6).

In case of an unlikely but very homogeneous situation, the normal-homogeneous approximation can be predominant generating a symmetric curve with mean μ and variance kμ.

*3.2.1 Analogy of the unlikely intensive to Log-normal distribution.*

It has been shown that when the experimental data are well represented by unlikely intensive natural distribution, they can be fitted also with the well known right skewed log-normal distribution. As mentioned above, the lack of phenomenological interpretation capability of the log-normal makes the intensive natural distribution the preferable one. It can be demonstrated that for not too dispersed data the two curves are very close to each other in a range from ½μ to 2μ.

Rewriting the unlikely intensive natural curve equation (5) and assuming $k = \sigma^2\mu$, consequently the resulting equation and the log-normal equation are

| Unlikely intensive natural | Log-normal |
| --- | --- |
| $f(x)=\frac{1}{\sigma\sqrt{2\pi\mu x}}\cdot e^{-\frac{1}{2\sigma^2}\left(\frac{x-\mu}{\sqrt{\mu x}}\right)^2}$ | $f(x)=\frac{1}{\sigma x\sqrt{2\pi}}\cdot e^{-\frac{1}{2\sigma^2}(\ln x-\ln\mu)^2}$ |

In this form it is postulated that the scale parameter σ and the position parameter μ have the same value. For the unlikely intensive natural distribution μ is the harmonic mean while for the log-normal one μ is the geometric mean.

Now we have found this useful approximation

$$\frac{x-\mu}{\sqrt{\mu x}} \approx ln\frac{x}{\mu}$$

then we demonstrate it considering the Taylor's series expansion of

$$\frac{1}{\sqrt{x}} = \frac{1}{\sqrt{\mu}} - \frac{1!!}{2!!}\frac{1}{\sqrt{\mu}}\left(\frac{x-\mu}{\mu}\right) + \frac{3!!}{4!!}\frac{1}{\sqrt{\mu}}\left(\frac{x-\mu}{\mu}\right)^2 - \frac{5!!}{6!!}\frac{1}{\sqrt{\mu}}\left(\frac{x-\mu}{\mu}\right)^3 + \frac{7!!}{8!!}\frac{1}{\sqrt{\mu}}\left(\frac{x-\mu}{\mu}\right)^4 - \dots$$

hence the left term can be expressed by this summation

$$\frac{x-\mu}{\sqrt{\mu x}} = \frac{x-\mu}{\sqrt{\mu}}\left[\frac{1}{\sqrt{\mu}} - \frac{1!!}{2!!}\frac{1}{\sqrt{\mu}}\left(\frac{x-\mu}{\mu}\right) + \frac{3!!}{4!!}\frac{1}{\sqrt{\mu}}\left(\frac{x-\mu}{\mu}\right)^2 - \frac{5!!}{6!!}\frac{1}{\sqrt{\mu}}\left(\frac{x-\mu}{\mu}\right)^3 + \frac{7!!}{8!!}\frac{1}{\sqrt{\mu}}\left(\frac{x-\mu}{\mu}\right)^4 - \dots\right]$$

which becomes

$$\frac{x-\mu}{\sqrt{\mu x}} = \frac{x-\mu}{\mu} - \frac{1!!}{2!!}\left(\frac{x-\mu}{\mu}\right)^2 + \frac{3!!}{4!!}\left(\frac{x-\mu}{\mu}\right)^3 - \frac{5!!}{6!!}\left(\frac{x-\mu}{\mu}\right)^4 + \frac{7!!}{8!!}\left(\frac{x-\mu}{\mu}\right)^5 - \dots$$

While the right term can be approximated using a Taylor expansion around μ as

$$ln\frac{x}{\mu} = \frac{x-\mu}{\mu} - \frac{1}{2}\left(\frac{x-\mu}{\mu}\right)^2 + \frac{1}{3}\left(\frac{x-\mu}{\mu}\right)^3 - \frac{1}{4}\left(\frac{x-\mu}{\mu}\right)^4 + \frac{1}{5}\left(\frac{x-\mu}{\mu}\right)^5 - \dots$$

Comparing the two sommatories, it is intuitive that the first two terms, *i.e.* the first and second order, coincide. Also the higher orders are similar. The **Figure 3** shows this similarity.

For values of x ranging between ½μ and 2μ, *i.e.* for x/μ in the range [0.5;2], $\frac{x-\mu}{\sqrt{\mu x}}$ and $ln\frac{x}{\mu}$ differ for less than 2% (2% on the extremes). In this range, the coefficient of determination is 98.4%; if the data are not too dispersed also $\sqrt{\mu x} \approx x$ and so the two frequency functions are similar. This closeness is mainly driven by the parameters of dispersion σ therefore, in case σ becomes small enough, the two curves become almost identical. For dispersed data the log-normal distribution is not a good approximation of the unlikely intensive one.

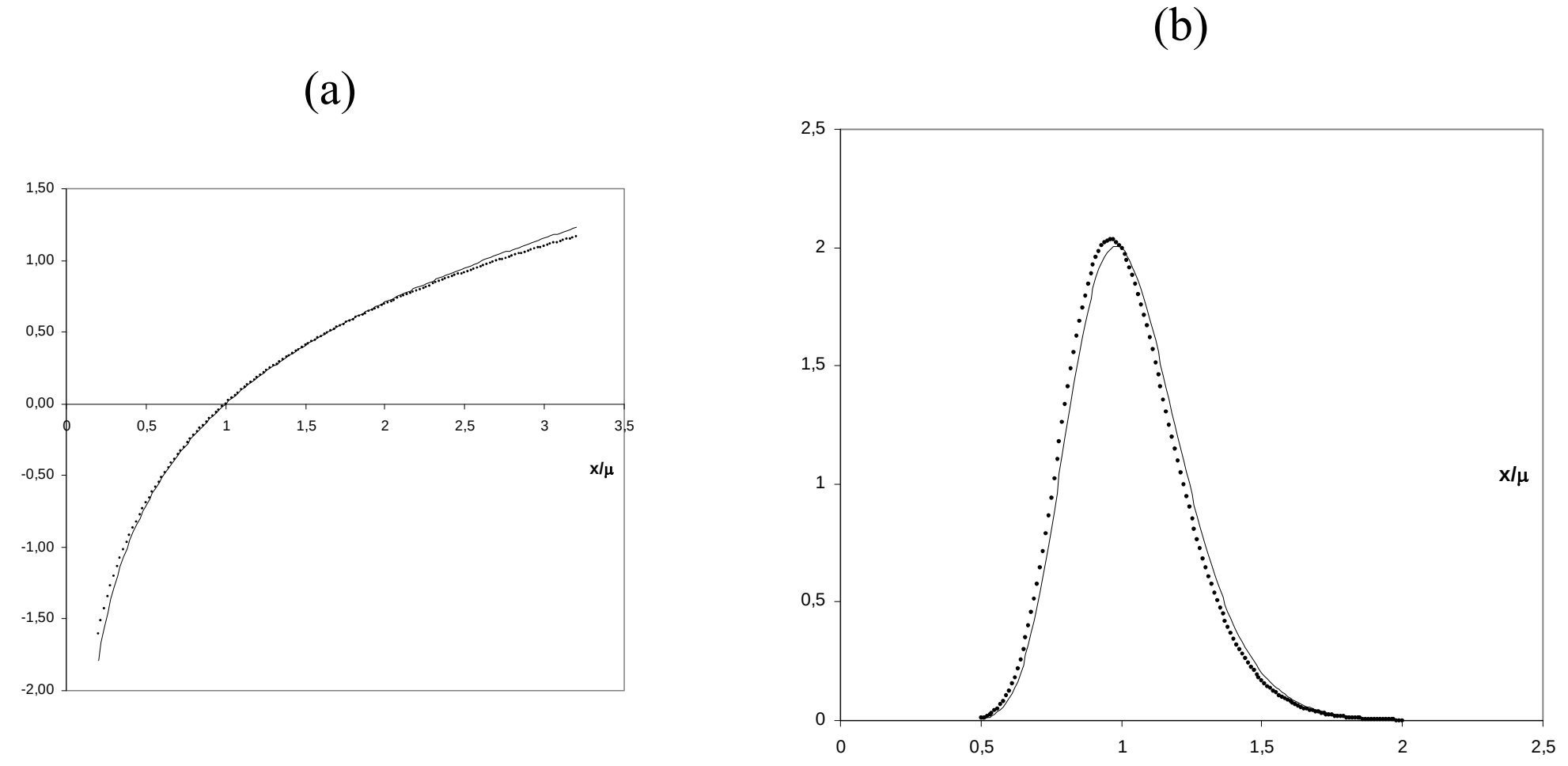

**Figure 3**. The similarity between log-normal and unlikely intensive natural functions: (a) plot of $\frac{x-\mu}{\sqrt{\mu x}}$ (—) and of $ln\frac{x}{\mu}$ (····) for various x/μ values; (b) plot of the unlikely intensive natural and lognormal frequency functions for σ = 0.2 and for various x/μ values.

*3.3 The likely intensive: a left skewed model.*

The likely approximation does not have any well-known correspondent discrete distribution. In this case, a left tailed frequency curve is less frequent to encounter; when the probability p tends to 1, the χ values having not zero probability are χ ≈ 1. This is like considering the generic equation (3) where x ≈ u. By this approximation and re-labeling k, the product ku, the frequency function (3) is written as

$$f(x)=\frac{1}{\sqrt{2\pi k(u-x)}}\cdot e^{-\frac{(\mu-x)^2}{2k(u-x)}} \qquad (6)$$

which is the mirrored representation of the unlikely approximation curve (5) symmetric respect to p = 0.5 (see **Figure 1**). In fact, substituting the variable x with y= u – x, the same equation (5), which estimates u – μ instead of μ (as overall mean value) and the same dispersion parameter k, is obtained. Accordingly, the same PROOFS are valid for both equation (5) and the mirrored likely intensive natural equation (6). As a result we can calculate the batch mean value μ and the dispersion value k

$$E\left(\frac{1}{u-x}\right)=\frac{1}{u-\mu} \qquad\qquad k=\mu-E(x)$$

Due to the left tail, the mean (between samples) is every time an underestimation of the overall mean value μ. The same conclusions for the unlikely intensive natural distribution can be drawn.

The variance is $\mathrm{var}(x)=k(2k+u-\mu)$ or $\mathrm{var}(x)=k^2+k(u-E(x))$.

The function (6) is not related to any discrete distribution, except for the binomial p when p tends to 1. When there is a likely event but very homogeneous, the normal-homogeneous approximation can be combined with the likely approximation generating a symmetric curve with mean (u-μ) and variance k(u-μ).

## 4. APPLICATIONS AND EXAMPLES

Illustrative applications, derived from chemical and geotechnical available datasets, are given (**Figures 4-8)**. We compared the proposed intensive natural distribution to other widely-used frequency functions to test its ability in fitting the histograms and the probability charts obtained from various intensive variables.

The parameters used for data reproduction are calculated on the same data using the methods shown in **Table 1**.

Alternatively, when the presence of a bimodal histogram or spurious data is detected, a better estimation of parameters is performed by best-fitting distribution over the histogram.

**Table 1.** List of formulas for parameterization of the frequency functions, with x the dataset array, E(x) the dataset mean and $\boldsymbol{\sigma}(\mathbf{x})$ the dataset standard deviation.

| Function type | parameters |
|---|---|
| unlikely intensive natural<br>log-normal | μ and k see par. 3.2<br>$\mu = e^{E(ln\ x)}$ and $\sigma = e^{\sigma(ln\ x)}$ |
| likely intensive natural<br>"mirrored" log-normal | μ and k see par. 3.3<br>$\mu = e^{E(ln\ x)}$ and $\sigma = e^{\sigma(ln\ x)}$, |
| homogeneous unlikely<br>normal | μ = E(x) and k = σ(x)/μ see par. 3.3<br>μ = E(x) and σ = σ(x) |

The mirrored log-normal distribution is simply a log-normal function mirrored obtained substituting x with u – y, to obtain a left skewed frequency function having the following equation

$$f(x) = \frac{1}{\sigma(u - x)\sqrt{2\pi}} \cdot e^{-\frac{1}{2\sigma^2}(ln(u-x) - ln(u-\mu))^2}.$$

All histograms are overlapped with the theoretical frequency functions. PP plots are obtained, from the growing ordered set of N data, plotting for each *datum* "*i*" the experimental probability calculated as $P_{i\,exp} = \frac{i - 0.5}{N}$ *versus* the calculated probability from the theoretical cumulative frequency function. PP plots are a quick tool to show the goodness of fit and the confidence level achieved.

**Figure 4** shows histogram and probability chart fitting of the results of 89 samples of about 1 kg each of agricultural soil collected for ground value determination. In intensive narural distribution inference, lead is randomly diffused and heterogeneous, so some areas are more enriched than others. In log-normal distribution inference, this heterogeneity cannot be determined.

<table>
<tr><th>Unlikely intensive natural</th><th>Log-normal</th></tr>
<tr><td>

(a)

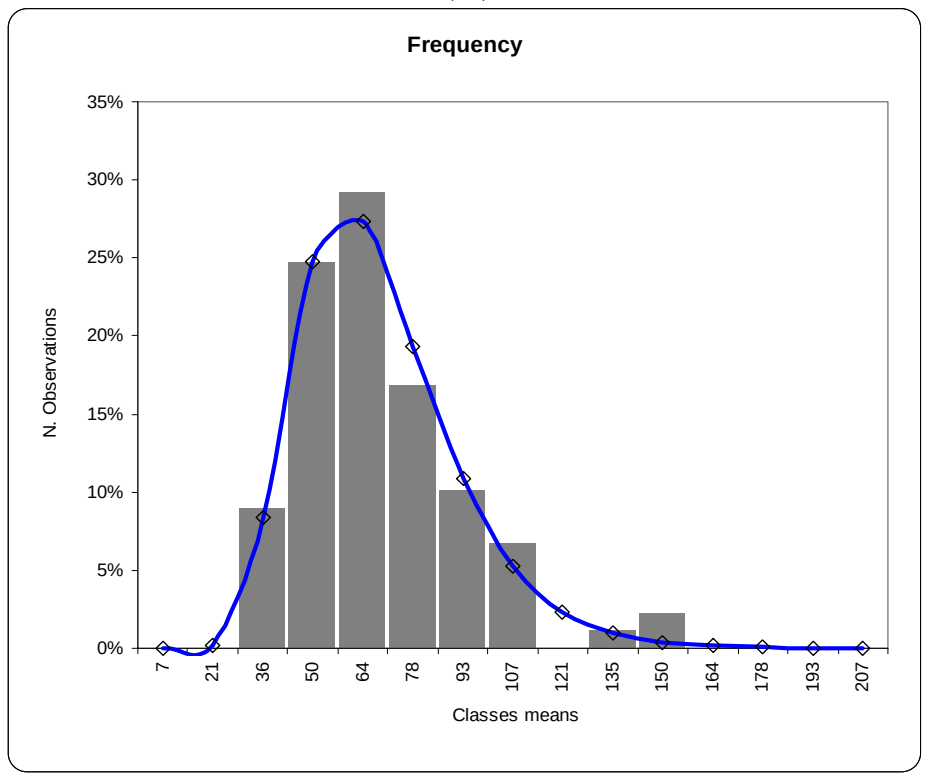

</td><td>

(b)

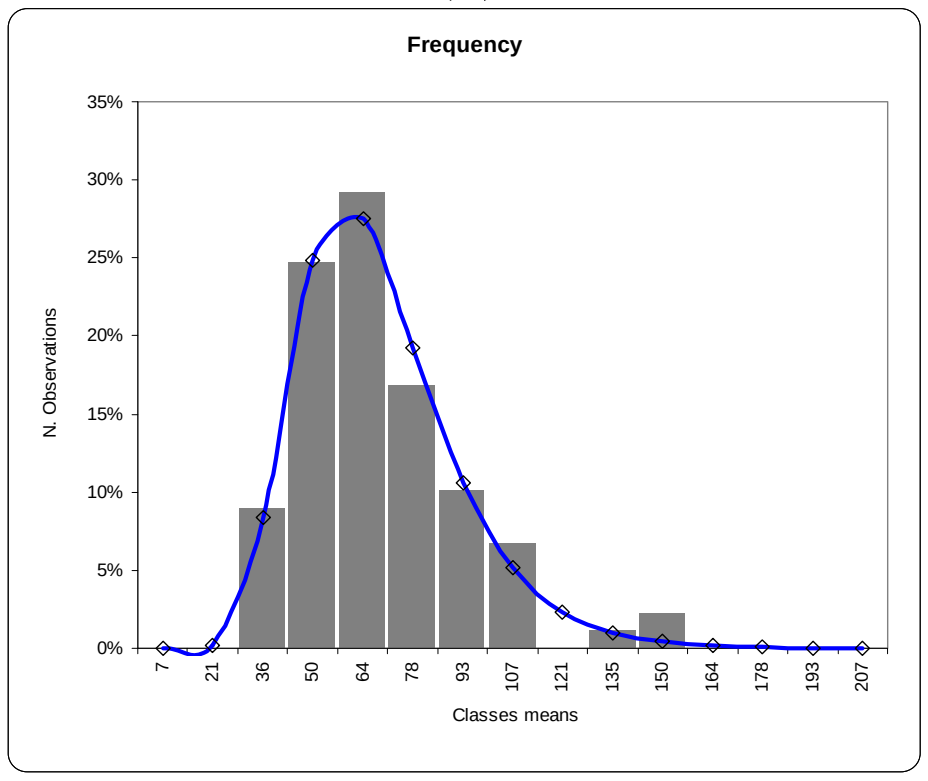

</td></tr>
<tr><td>

(c)

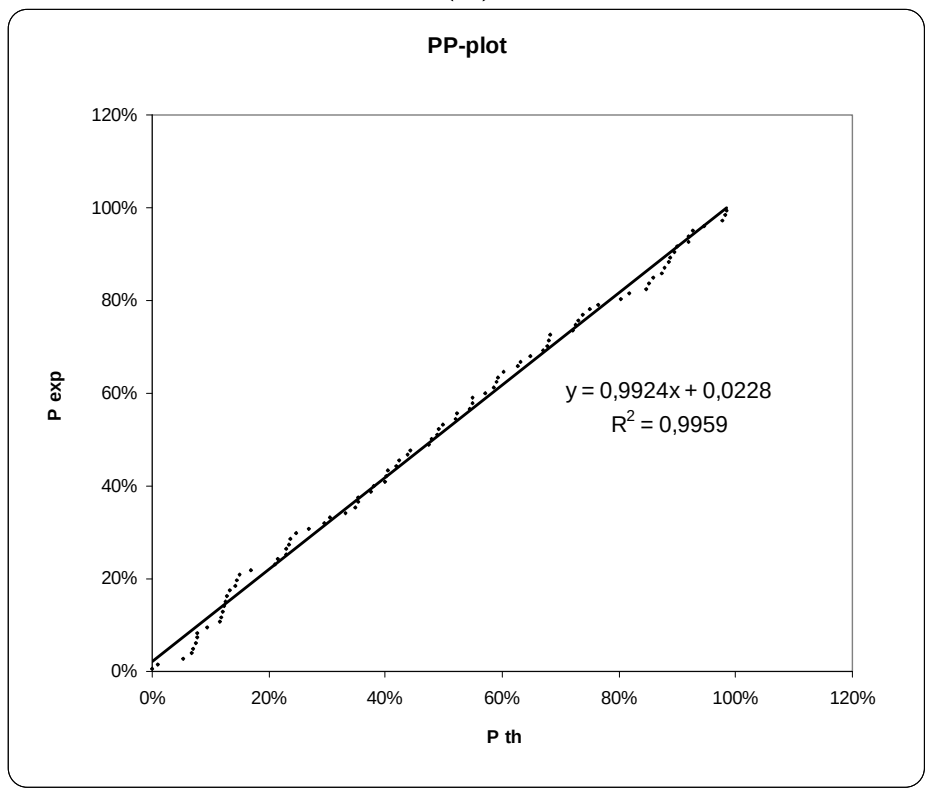

</td><td>

(d)

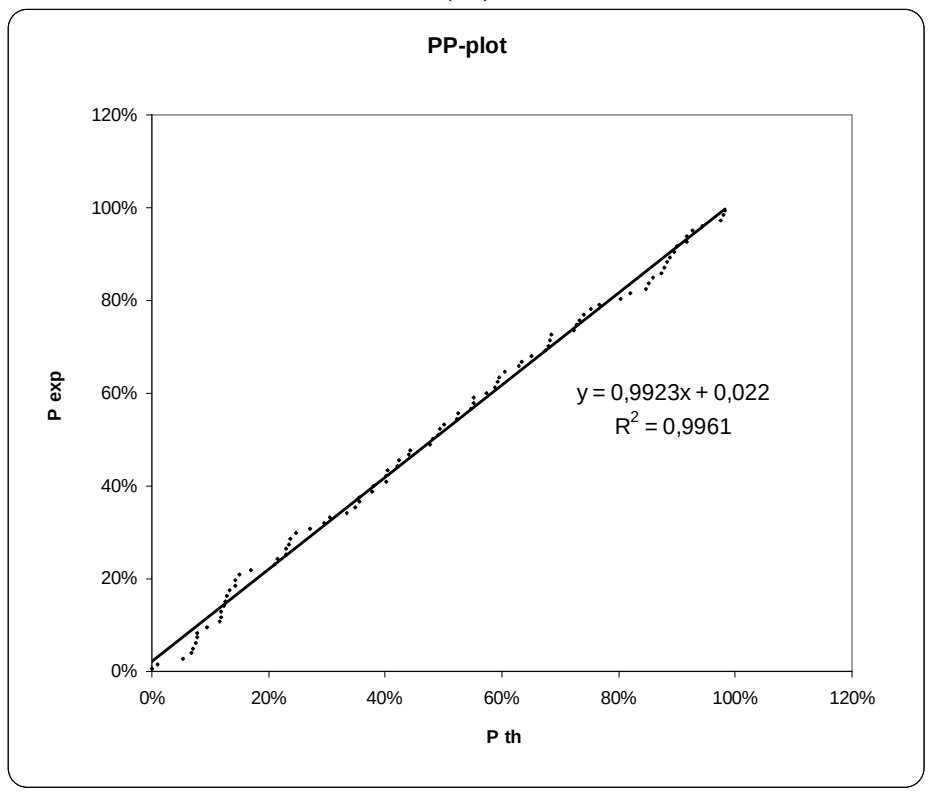

</td></tr>
<tr><td>$\mu = 63.2$<br>$k = 6.70$</td><td>$\mu = 64.1$<br>$\sigma = 0.311$</td></tr>
</table>

**Figure 4**. The four plots show the goodness of fit of lead experimental data measurements of in agricultural soils (histogram fitting with (a) unlikely intensive natural model and with (b) log-normal model and probability chart (PP plot) for (c) unlikely intensive natural probability and with (d) log-normal probability).

**Figure 5** shows histogram and probability chart fitting of the results of 169 samples of about 10 mL each of underground water collected for ground value determination. Manganese was measured and with intensive natural distribution inference, it was found randomly diffused and heterogeneous in space and or time. Using the log-normal distribution inference, no conclusions can be drawn. The parameters are better reproduced by the best- histogram fitting than from the maximum likehood estimation. In this case the maximum likehood estimated parameters are used to show the over-dispersion due to the presence of some spurious data.

<table>
<tr><th>Unlikely intensive natural</th><th>Log-normal</th></tr>
<tr><td>(a)<br>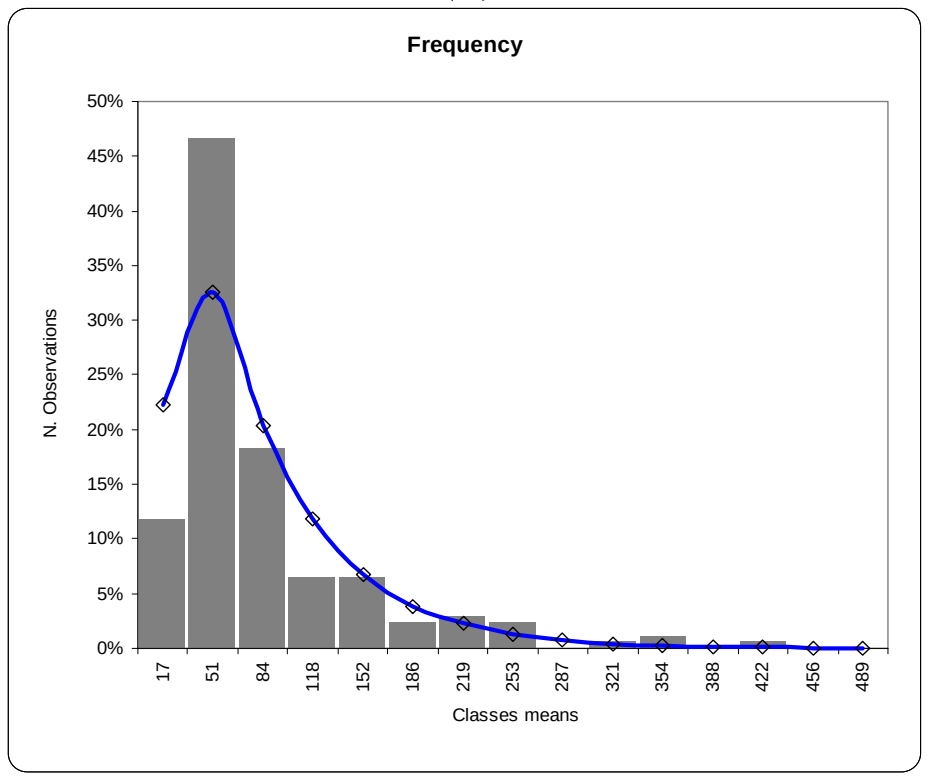
</td><td>(b)<br>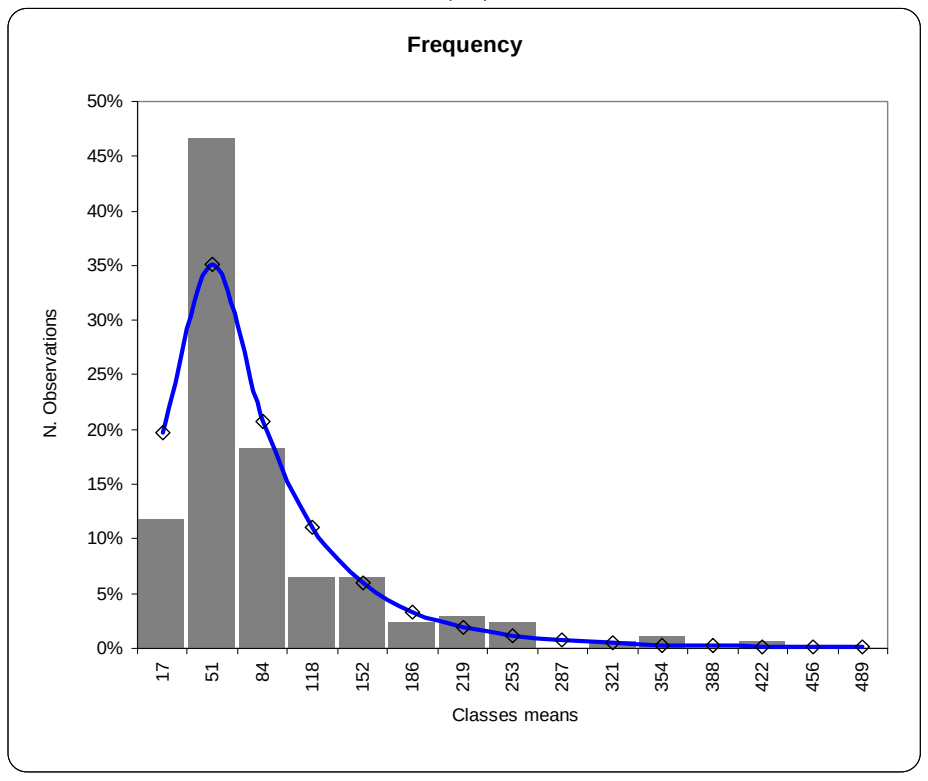
</td></tr>
<tr><td>(c)<br>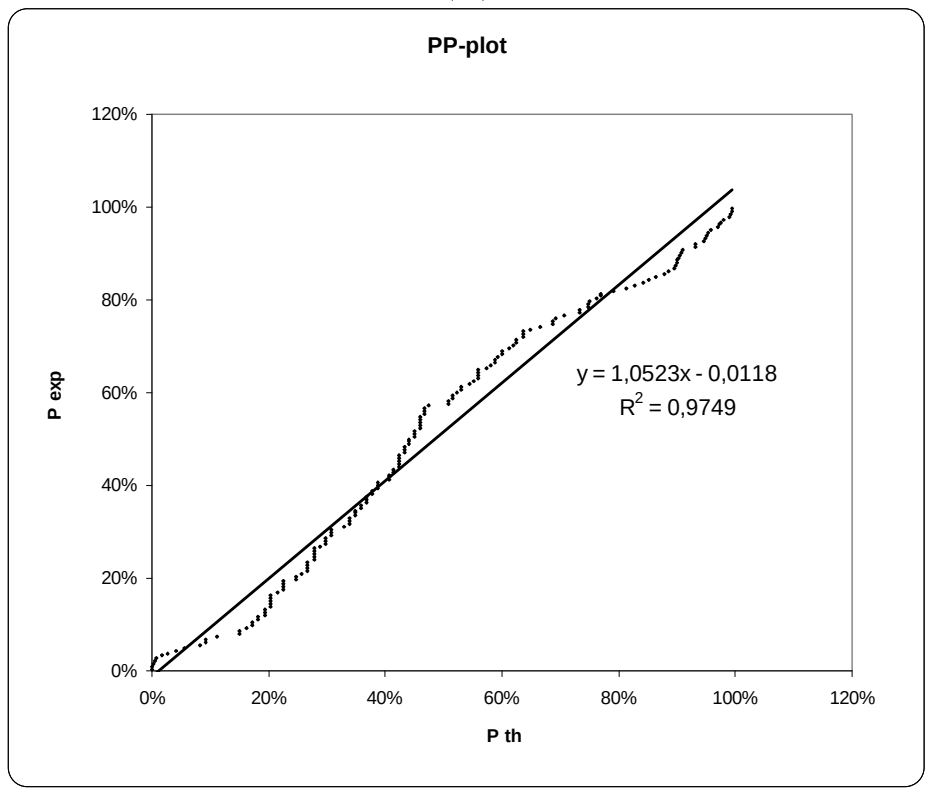
</td><td>(d)<br>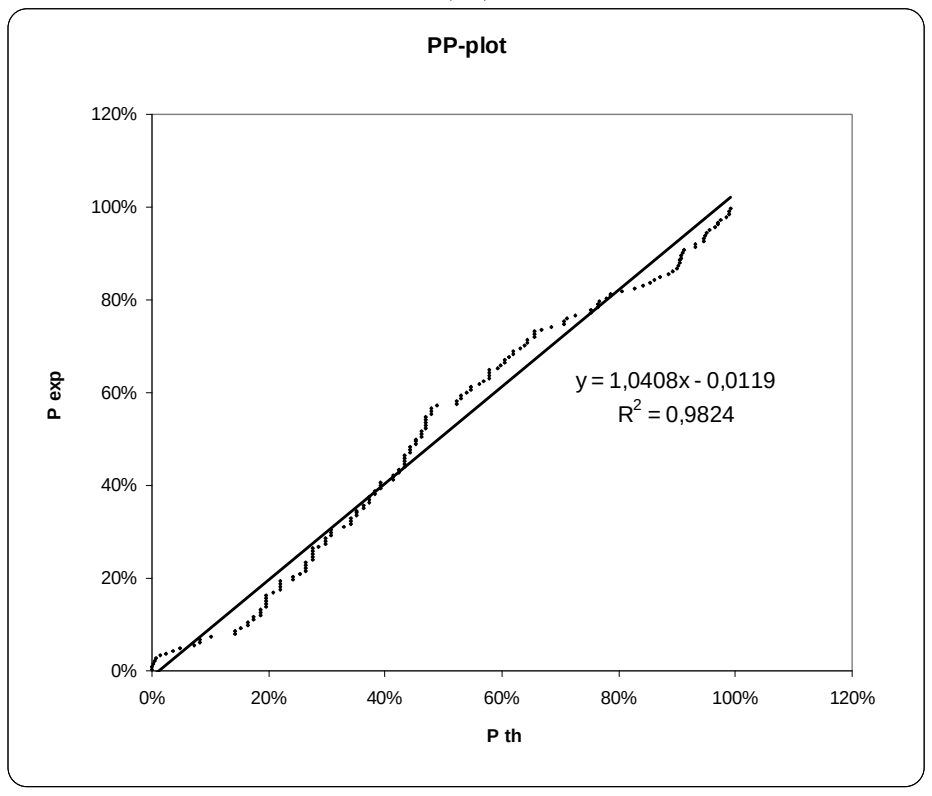
</td></tr>
<tr><td>$\mu = 49.5$<br>$k = 33.8$</td><td>$\mu = 64.1$<br>$\sigma = 2.05$</td></tr>
</table>

**Figure 5**. The four plots show the goodness of fit of manganese experimental measurements data in underground waters (histogram fitting with (a) unlikely intensive natural model and with (b) log-normal model and probability chart (PP plot) for (c) unlikely intensive natural probability and with (d) log-normal probability).

**Figure 6** shows histogram and probability chart fitting of the results of approximately one year of PM10 monitoring (342 data, one per day) for atmospheric air quality, collected by a stationary automated measuring system. In intensive natural distribution inference, PM10 is randomly present and heterogeneous in time. Using log-normal distribution inference, these conclusions cannot be determined.

<table>
<tr><td>Unlikely intensive natural</td><td>Log-normal</td></tr>
<tr><td>(a)<br>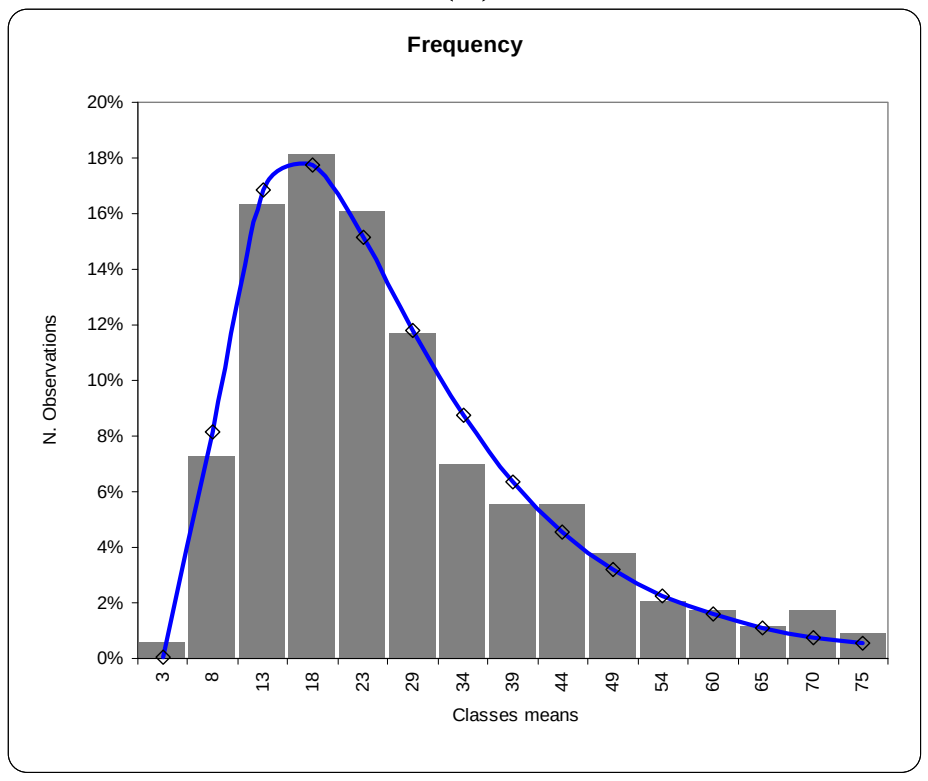
</td><td>(b)<br>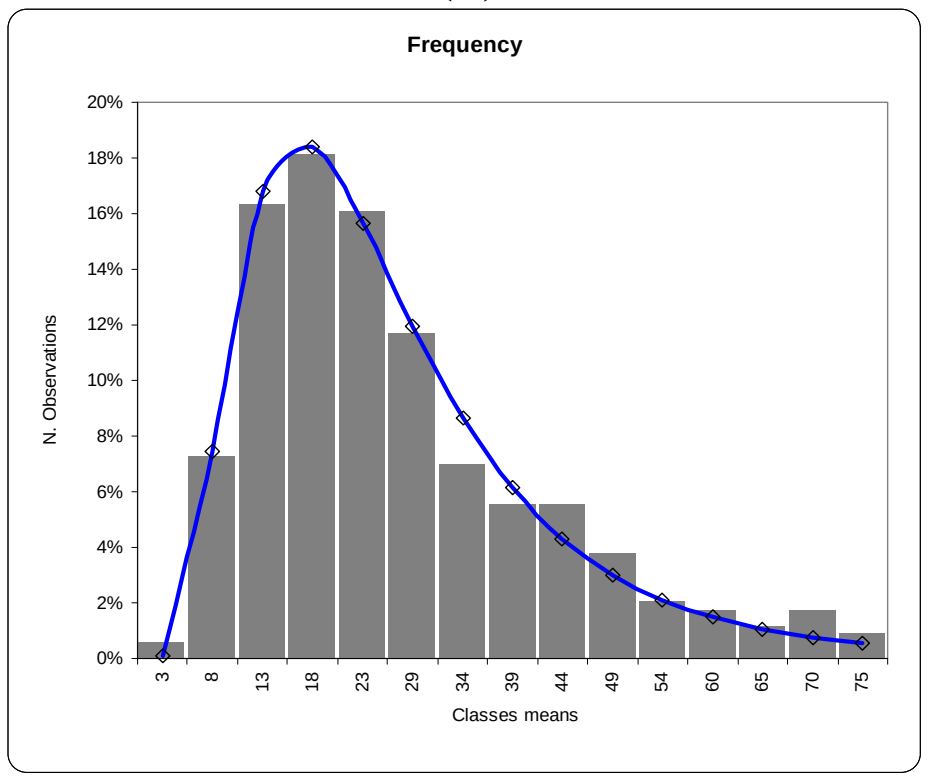
</td></tr>
<tr><td>(c)<br>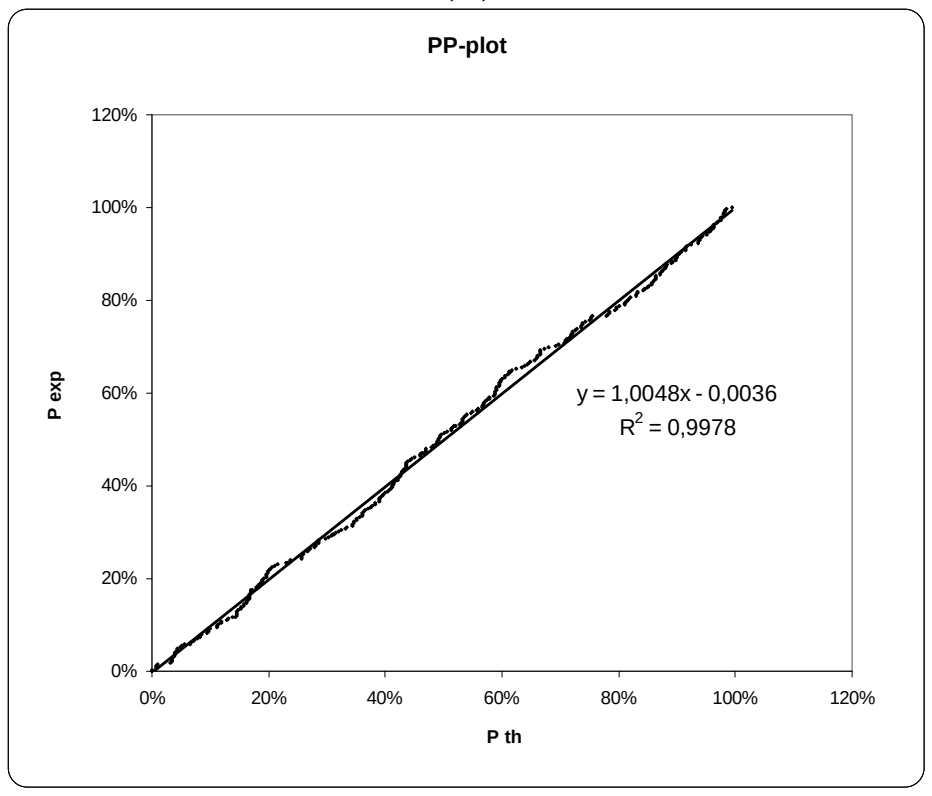
</td><td>(d)<br>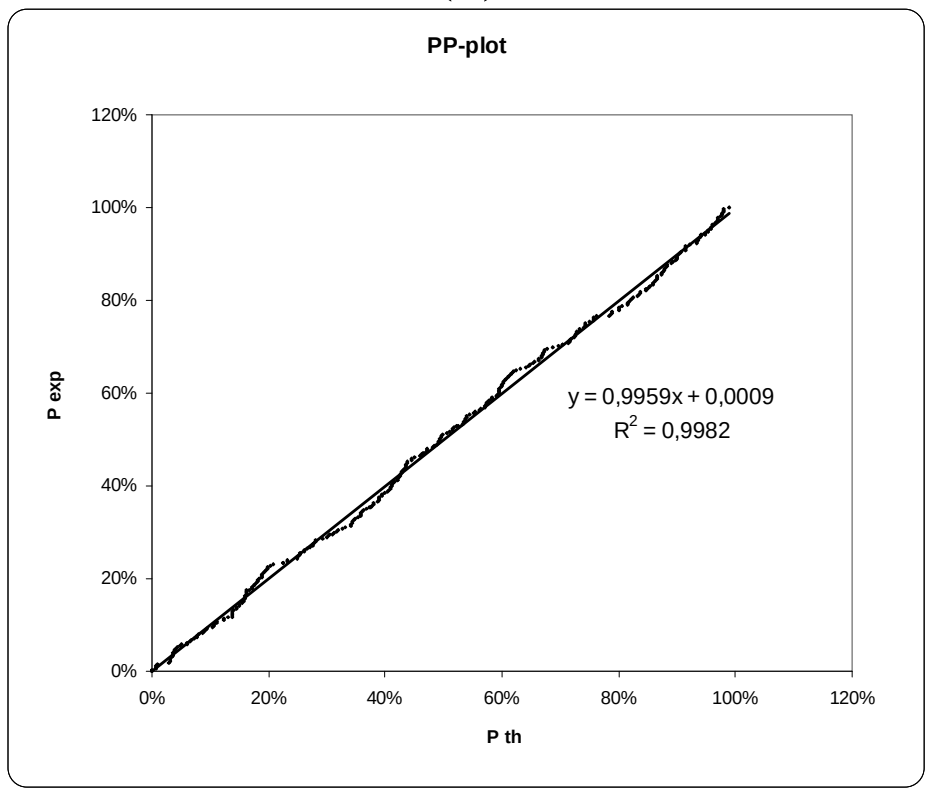
</td></tr>
<tr><td>$\mu = 19.6$<br>$k = 7.63$</td><td>$\mu = 23.1$<br>$\sigma = 1.76$</td></tr>
</table>

**Figure 6**. The four plots show the goodness of fit of PM10 experimental data in atmospheric air (histogram fitting with (a) unlikely intensive natural model and with (b) log-normal model and probability chart (PP plot) for (c) unlikely intensive natural probability and with (d) log-normal probability).

**Figure 7** shows histogram and probability chart fitting of the results of 180 samples of 1 Kg lagoon sediments collected for ground value determination of the lagoon. Chromium was measured and found distributed random but homogeneously present in all the area.

| Unlikely homogeneous intensive natural | Normal (Gaussian) |
|---|---|
| (a) 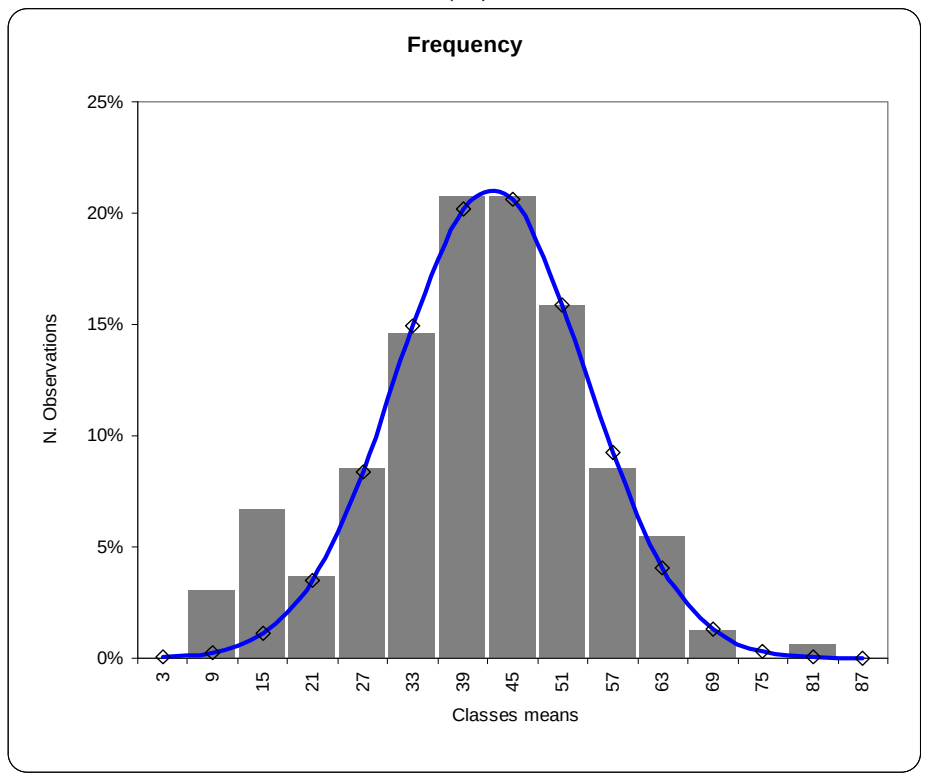 | (b) 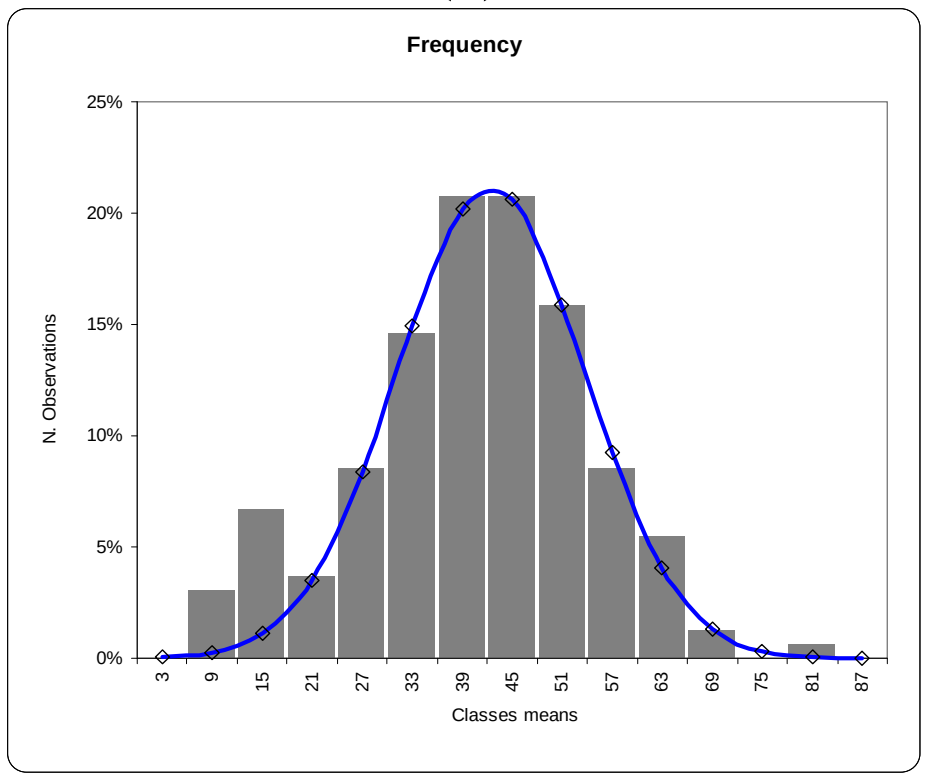 |
| (c) 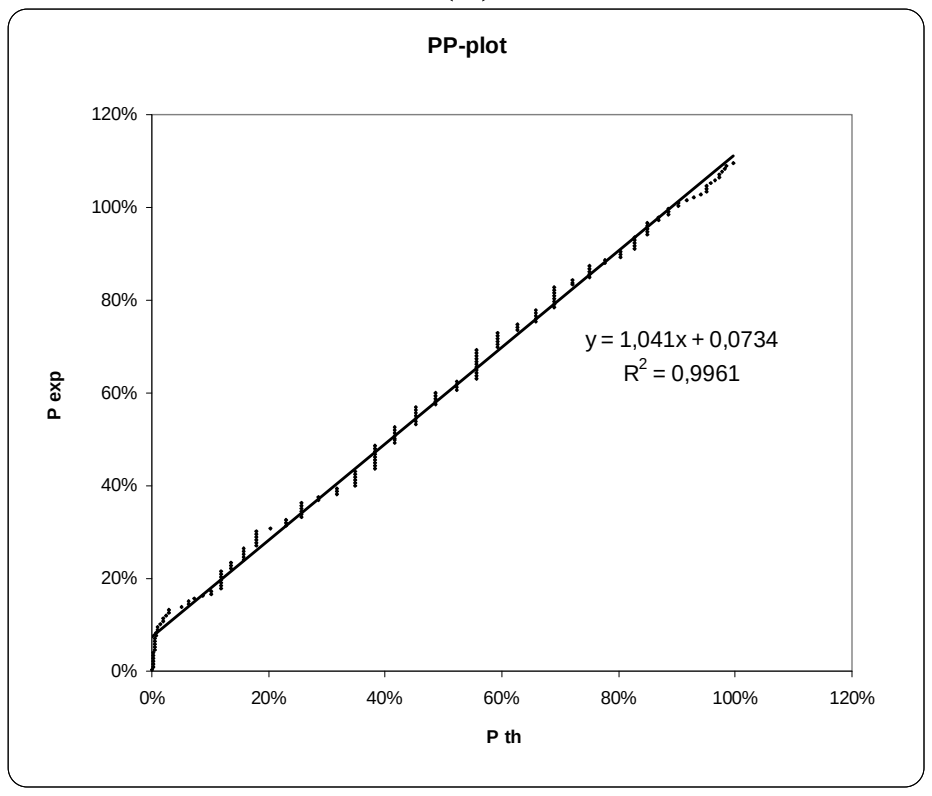 | (d) 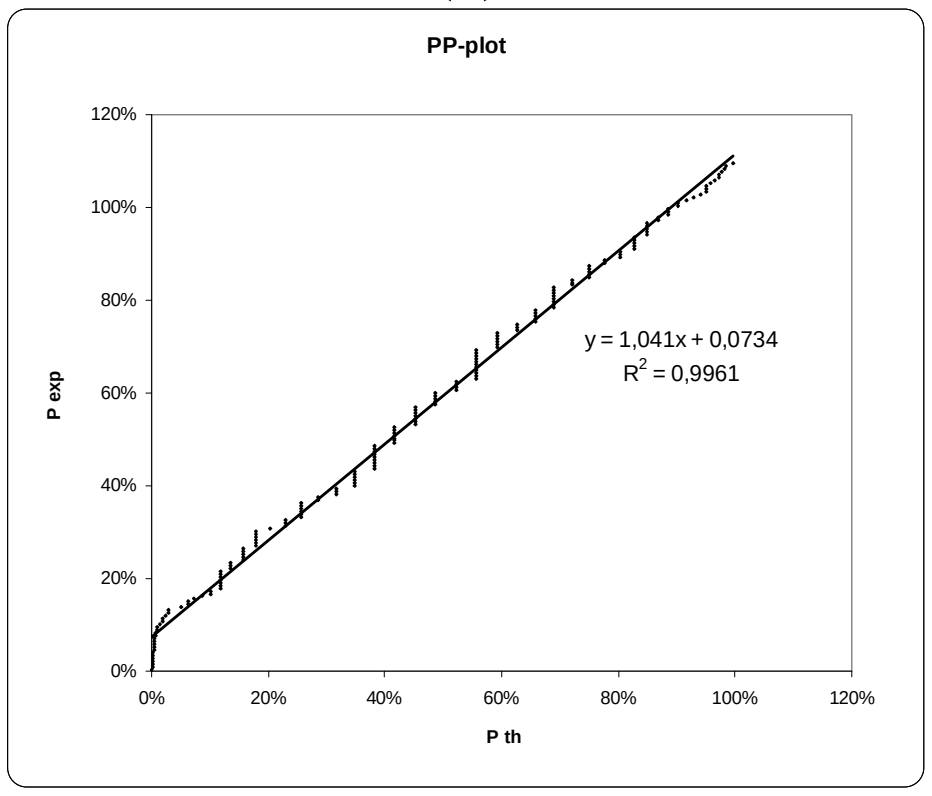 |
| $\mu = 42.4$<br>$k = 3.02$ | $\mu = 42.4$<br>$\sigma = 11.3$ |

**Figure 7**. The four plots show the goodness of fit of chromium experimental measurements data from lagoon sediments (histogram fitting with (a) unlikely intensive natural model and with (b) log-normal model and probability chart (PP plot) for (c) unlikely intensive natural probability and with (d) log-normal probability).

The parameters are estimated by best fitting the histogram with the theoretical curve to minimize the effect of a small second population which causes over dispersion. Being the homogeneity related to the entropy and the entropy to the time passed, it can be concluded that the presence of Cr in sediment is historical. In fact it was all the time necessary to diffuse. Using the log-normal distribution inference, a mean value of Cr randomly present in the area is determined.

**Figure 8** shows histogram and probability chart fitting of the results of the geotechnical analysis of silt percentage in about 4000 samples of sediments determined by laser scattering technique.

<table>
<tr><td>Likely intensive natural</td><td>“mirrored” Log-normal</td></tr>
<tr><td>
(a)
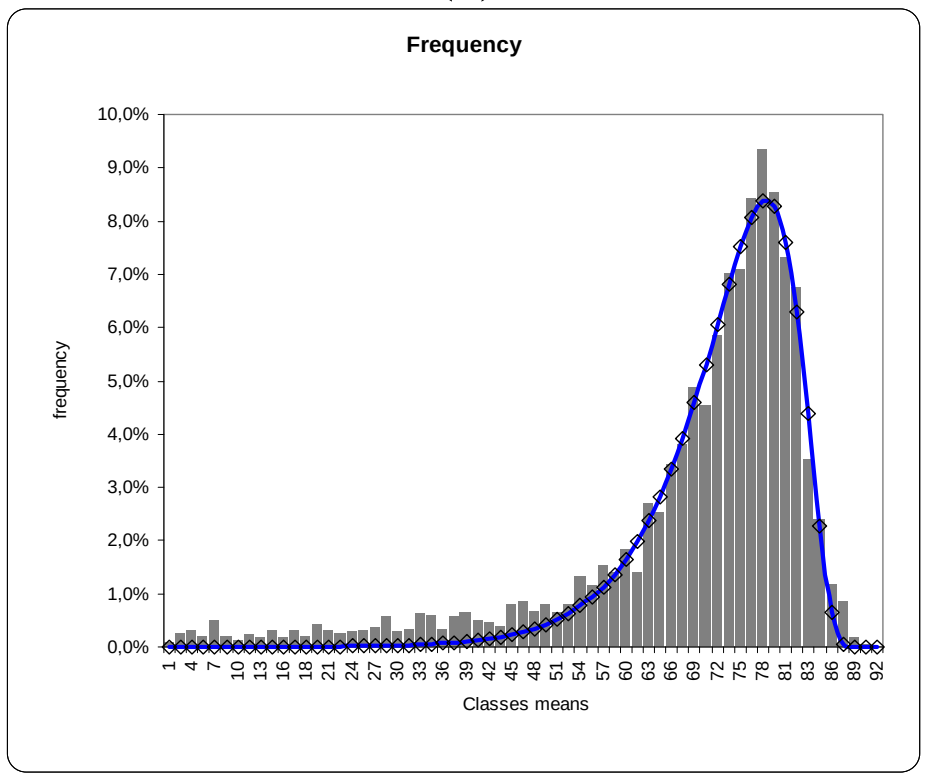

</td><td>
(b)
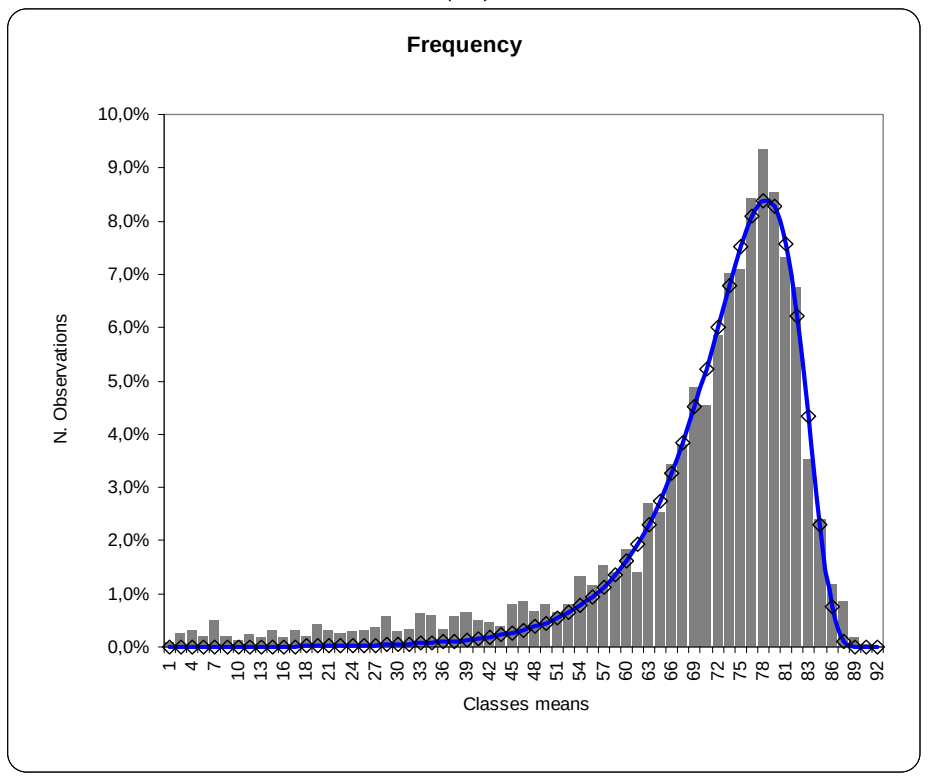

</td></tr>
<tr><td>
(c)
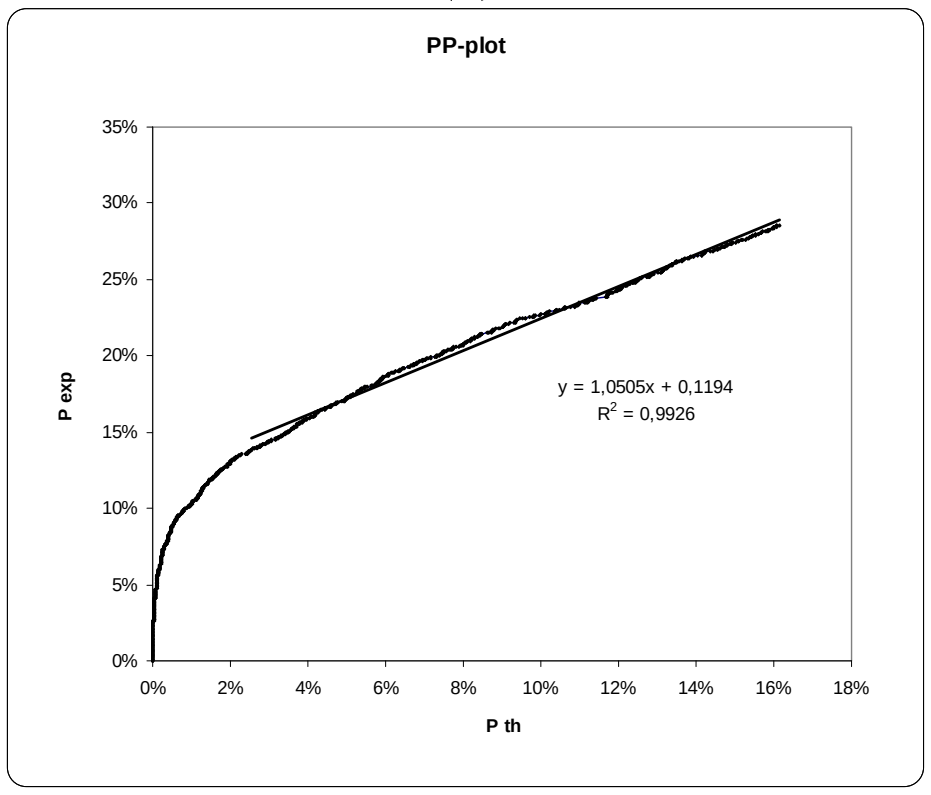

</td><td>
(d)
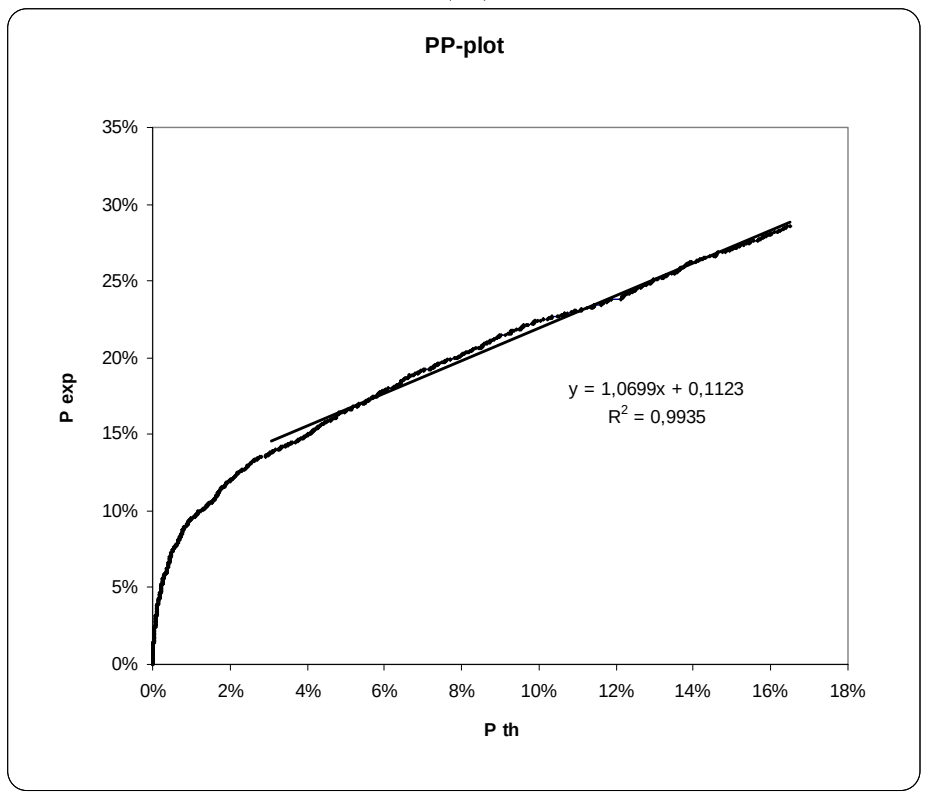

</td></tr>
</table>

**Figure 8.** The four plots show the goodness of fit of % silt experimental measurements data of in sediments (histogram fitting with (a) unlikely intensive natural model and with (b) log-normal model and probability chart (PP plot) for (c) unlikely intensive natural probability and with (d) log-normal probability. The initial trend in PP plots is due to the presence of aggregates of clay or other particle size overestimated which is typical of the laser scattering method used).

The parameters are estimated by best fitting the histogram with the theoretical curve to minimize the effect of a series of spurious data present at low %. The percentage of % of silt was randomly and heterogeneously present in sediments.

**5. DISCUSSION**

Gaussian law has been used since the time of Laplace as the continuous approximation of the binomial distribution although it was evident that the asymmetry capability was lost by the Gaussian model approximation. Nowadays, this limitation is becoming increasingly apparent and detrimental due to the huge amount of data and to the availability of high resolution instrumentations. The explicit heterogeneity of the measurements becomes increasingly manifest as the ability to detect chemicals at lower concentrations than in the past progresses, and so does the need for robust asymmetrical distributions. The ability of detecting chemicals at progressively lower concentrations stresses more than in the past the heterogeneity of the measurements and therefore the need of robust asymmetrical distributions. The statisticians developed several classes of distributions capable of skewness but there were mostly artificially constructed and therefore not strictly linked to a solid inference of the causes of the phenomenon.

The proposed approach is the final development to the continuous field for independent random intensive variables, also applicable to the extremes probabilities which were neglected in the Gaussian model derivation.

In this paper, all the variability has been attributed exclusively to the between-sample variability, in which each sample is a sample or subsample of the batch to test, and its variability increases with the heterogeneity of the batch and decreases if the sample to analyze tends to the whole batch. In fact, if the sample is the entire batch, then the variability will be zero while in case the sample will be infinitesimal respect to the whole batch, the variability will be maximum.

The developed model can be easily applied, and allows to explicit the $\sigma$ of the Gaussian model explaining the origin of the Horwitz equation. Several other interesting explanations of well-known scientific issues in different scientific fields, based on the application of this model, are currently under study.

**ACKNOWLEDGEMENT**

We thank the ARPA FVG lab's chief director Dr. Stefano Pison and the director of the provincial lab of Udine Dr Anna Lutman for their support. We are grateful to all the technicians of the lab of Udine for the analytical data and a special thank to Elisa Piccoli and Denis Mazzilis for heavy metals analysis.

APPENDIX

*PROOF 1: The normalization condition of the unlikely intensive natural distribution*

Let allow demonstrating that $\int_0^{+\infty} f(x)\,dx = 1$.

$$\int_0^{+\infty} f(x)\,dx = \int_0^{+\infty} \frac{1}{\sqrt{2\pi kx}} \cdot e^{-\frac{(x-\mu)^2}{2kx}}\,dx$$

$$= \frac{1}{\sqrt{2\pi}} \cdot \int_0^{+\infty} \frac{1}{\sqrt{kx}} \cdot e^{-\frac{\mu^2}{2kx}} \cdot e^{\frac{\mu}{k}} \cdot e^{-\frac{x}{2k}}\,dx \qquad \text{calling } Z = \sqrt{\frac{x}{2k}}$$

$$= \frac{2\sqrt{2}}{\sqrt{2\pi}} \cdot e^{\frac{\mu}{k}} \cdot \int_0^{+\infty} e^{-\frac{\mu^2}{4kZ^2}} \cdot e^{-Z^2}\,dZ$$

$$= \frac{2}{\sqrt{\pi}} \cdot e^{2a} \cdot \int_0^{+\infty} e^{-\left(\frac{a}{Z}\right)^2} \cdot e^{-Z^2}\,dZ \qquad \text{where } a = \frac{\mu}{2k}$$

The solution of the last integral is known

$$\int_0^{+\infty} e^{-\left(\frac{a}{Z}\right)^2} \cdot e^{-Z^2}\,dZ = \frac{\sqrt{\pi}}{4}\left[e^{-2|a|}\left\{1-\operatorname{erf}\left(\frac{|a|}{Z}-Z\right)\right\} - e^{2|a|}\left\{1-\operatorname{erf}\left(\frac{|a|}{Z}+Z\right)\right\}\right]_0^{+\infty}$$

where $\operatorname{erf}(x) = \frac{2}{\sqrt{\pi}}\int_0^x e^{-t^2}dt$ is the error function so $\operatorname{erf}(-\infty) = -1$ and $\operatorname{erf}(+\infty) = 1$.

Considering that is a positive $\lim\limits_{Z\to 0^+}\left(\frac{a}{Z} \pm Z\right) = +\infty$ and $\lim\limits_{Z\to +\infty}\left(\frac{a}{Z} \pm Z\right) = \pm\infty$.

The integral is the equal to $\int_0^{+\infty} e^{-\left(\frac{a}{Z}\right)^2} \cdot e^{-Z^2}\,dZ = \frac{\sqrt{\pi}}{2}e^{-2a}$

Finally $\int_0^{+\infty} f(x)\,dx = \frac{2}{\sqrt{\pi}} \cdot e^{2a} \cdot \frac{\sqrt{\pi}}{2} \cdot e^{-2a} = 1$.

*PROOF 2: The mean batch value $\mu$ of the unlikely intensive natural distribution*

Solving the integral E(x) is really a challenging problem, but it can be found the solution for E(1/x). Let us consider that

$$E\left(\frac{x+\mu}{2x}\right)=\int_0^{+\infty}\frac{x+\mu}{2x}\frac{1}{\sqrt{2\pi kx}}\cdot e^{-\frac{(x-\mu)^2}{2kx}}\,dx$$

$$=\frac{1}{\sqrt{2\pi}}\cdot\int_{-\infty}^{+\infty}\frac{x+\mu}{2x}\cdot\frac{2x}{x+\mu}\cdot e^{-\frac{Z^2}{2}}\,dZ \qquad \text{where } Z=\frac{\mu-x}{\sqrt{kx}}$$

$$=\frac{1}{\sqrt{2\pi}}\cdot\int_{-\infty}^{+\infty}e^{-\frac{Z^2}{2}}\,dZ$$

$$=1$$

Then $E\left(\frac{\mu+x}{2x}\right)=1$ $\quad E\left(\frac{\mu}{x}\right)=1$ $\quad$ and finally $\quad E\left(\frac{1}{x}\right)=\frac{1}{\mu}$.

*PROOF 3: The dispersion factor k of the unlikely intensive natural distribution*

The expression of the k factor was theoretically found starting from the normalization condition (see PROOF 1)

$$1=\int_0^{+\infty}\frac{1}{\sqrt{2\pi kx}}\cdot e^{-\frac{(x-\mu)^2}{2kx}}\,dx=\int_0^{+\infty}1\cdot\frac{1}{\sqrt{2\pi kx}}\cdot e^{-\frac{(x-\mu)^2}{2kx}}\,dx$$

and using the rule of the partial integration $\int u'v=uv-\int u'v$ whit u′ = 1 and v the remaining part

$$=\frac{x}{\sqrt{2\pi kx}}\cdot e^{-\frac{(x-\mu)^2}{2kx}}\Bigg|_0^{+\infty}-\int_0^{+\infty}\frac{x}{\sqrt{2\pi k}}\cdot\left\{-\frac{1}{2x\sqrt{x}}\cdot e^{-\frac{(x-\mu)^2}{2kx}}-\frac{1}{\sqrt{x}}\cdot e^{-\frac{(x-\mu)^2}{2kx}}\cdot\frac{x^2-\mu^2}{2kx^2}\right\}dx$$

the first term is zero because

$$\lim_{x\to+\infty}\frac{x}{\sqrt{2\pi kx}}\cdot e^{-\frac{(x-\mu)^2}{2kx}}=\frac{1}{\sqrt{2\pi k}}\cdot\lim_{x\to+\infty}\frac{\sqrt{x}}{e^{\frac{x}{2k}}}=\frac{1}{\sqrt{2\pi k}}\cdot\lim_{x\to+\infty}\frac{1}{2\sqrt{x}\cdot\frac{1}{2k}e^{\frac{x}{2k}}}=\sqrt{\frac{k}{2\pi}}\cdot\frac{1}{\infty}=0$$

and $$\lim_{x\to0}\frac{x}{\sqrt{2\pi kx}}\cdot e^{-\frac{(x-\mu)^2}{2kx}}=\frac{1}{\sqrt{2\pi k}}\cdot\lim_{x\to+\infty}\frac{\sqrt{x}}{e^{\frac{\mu^2}{2kx}}}=\frac{1}{\sqrt{2\pi k}}\cdot\frac{0}{\infty}=0$$

The second term can be divided in two parts

$$1=\frac{1}{2}\int_0^{+\infty}\frac{1}{\sqrt{2\pi kx}}\cdot e^{-\frac{(x-\mu)^2}{2kx}}dx+\frac{1}{2}\int_0^{+\infty}\frac{x^2-\mu^2}{kx}\cdot\frac{1}{\sqrt{2\pi kx}}\cdot e^{-\frac{(x-\mu)^2}{2kx}}dx$$

the first one is again the normalization integral and equals to 1; the second integral is the expectation value of a function, i.e.

$$E\left(\frac{x^2-\mu^2}{kx}\right)=\int_0^{+\infty}\frac{x^2-\mu^2}{kx}\cdot\frac{1}{\sqrt{2\pi kx}}\cdot e^{-\frac{(x-\mu)^2}{2kx}}dx\,.$$

So the equation became

$$1=\frac{1}{2}+\frac{1}{2}E\left(\frac{x^2-\mu^2}{kx}\right)\qquad E\left(\frac{x^2-\mu^2}{kx}\right)=1\qquad E\left(x-\frac{\mu^2}{x}\right)=k$$

and then $$k=E(x)-\mu^2E\left(\frac{1}{x}\right)$$

Finally using the expression of the expectation value (see PROOF 2)

$$k=E(x)-\mu$$

and $$E(x)=k+\mu\,.$$

*PROOF 4: The variance of the unlikely intensive natural distribution*

The expression of var(x) was theoretically found starting from the expectation value definition (see PROOF 2)

$$E(x)=\int_0^{+\infty} x\frac{1}{\sqrt{2\pi kx}}\cdot e^{-\frac{(x-\mu)^2}{2kx}}\,dx$$

and using the rule of the partial integration $\int u'v = uv - \int u'v$ whit $u'$ = x and v the remaining part

$$=\frac{x^2}{2}\cdot\frac{1}{\sqrt{2\pi kx}}\cdot e^{-\frac{(x-\mu)^2}{2kx}}\Bigg|_0^{+\infty} - \int_0^{+\infty}\frac{x^2}{2}\cdot\frac{1}{\sqrt{2\pi k}}\cdot\left(-\frac{1}{2x\sqrt{x}}\cdot e^{-\frac{(x-\mu)^2}{2kx}} - \frac{1}{\sqrt{x}}\cdot e^{-\frac{(x-\mu)^2}{2kx}}\cdot\frac{x^2-\mu^2}{2kx^2}\right)dx$$

the first term is zero because (see also PROOF 2)

$$\lim_{x\to+\infty}\frac{x^2}{2\sqrt{2\pi kx}}\cdot e^{-\frac{(x-\mu)^2}{2kx}} = \frac{1}{2\sqrt{2\pi k}}\cdot\lim_{x\to+\infty}\frac{x^{\frac{3}{2}}}{e^{\frac{x}{2k}}} = \frac{-3}{4\sqrt{2\pi k}}\cdot\lim_{x\to+\infty}\frac{\sqrt{x}}{\frac{1}{2k}e^{\frac{x}{2k}}} = 0$$

and
$$\lim_{x\to 0}\frac{x^2}{\sqrt{2\pi kx}}\cdot e^{-\frac{(x-\mu)^2}{2kx}} = \frac{1}{\sqrt{2\pi k}}\cdot\lim_{x\to+\infty}\frac{x^{\frac{3}{2}}}{e^{\frac{\mu^2}{2kx}}} = \frac{1}{\sqrt{2\pi k}}\cdot\frac{0}{\infty} = 0$$

The second term can be divided in two parts

$$E(x)=\frac{1}{4}\int_0^{+\infty} x\frac{1}{\sqrt{2\pi kx}}\cdot e^{-\frac{(x-\mu)^2}{2kx}}\,dx + \frac{1}{4k}\int_0^{+\infty}\left(x^2-\mu^2\right)\frac{1}{\sqrt{2\pi kx}}\cdot e^{-\frac{(x-\mu)^2}{2kx}}\,dx$$

the first one is again the E(x) and the second integral is the expectation value of a function, *i.e.*

$$E\left(x^2-\mu^2\right)=\int_0^{+\infty}\left(x^2-\mu^2\right)\frac{1}{\sqrt{2\pi kx}}\cdot e^{-\frac{(x-\mu)^2}{2kx}}\,dx\ .$$

So the equation became

$$E(x)=\frac{1}{4}E(x)+\frac{1}{4}E\left(x^2-\mu^2\right) \qquad 3kE(x)=E\left(x^2\right)-\mu^2$$

$E(x^2) = 3kE(x) + \mu^2 = 3k^2 + 3k\mu + \mu^2$ (see PROOF 2).

Then by definition of variance

$$\mathrm{var}(x) = E(x^2) - E(x)^2 = 3k^2 + 3k\mu + \mu^2 - (k+\mu)^2$$

$\mathrm{var}(x) = 2k^2 + k\mu$ or $\mathrm{var}(x) = k^2 + kE(x)$.

*PROOF 5: The moment generating and the characteristic function for the unlikely intensive natural distribution*

The moment generating function by definition is obtained solving the following integral

$$M_x(t) = E(e^{tx}) = \int_0^{+\infty} e^{tx} \frac{1}{\sqrt{2\pi kx}} e^{-\frac{(x-\mu)^2}{2kx}} dx$$

the last integral equals to

$M_x(t) = \int_0^{+\infty} \frac{1}{\sqrt{2\pi kx}} e^{-\left[\frac{(x-\mu)^2}{2kx} - tx\right]} dx$ where the argument of the exponential is the critical part and it can be developed as follows

| | |
|---|---|
| $\frac{(x-\mu)^2}{2kx} - tx = \frac{x^2 - 2x\mu + \mu^2 - 2ktx^2}{2kx}$ <br> $= \frac{x^2(1-2kt) - 2x\mu + \mu^2 - 2x\mu\sqrt{1-2kt} + 2x\mu\sqrt{1-2kt}}{2kx}$ <br> $= \frac{\left(x\sqrt{1-2kt} - \mu\right)^2 - 2x\mu\left(1 - \sqrt{1-2kt}\right)}{2kx}$ <br> $= \frac{\left(x\sqrt{1-2kt} - \mu\right)^2}{2kx} - \frac{\mu}{k}\left(1 - \sqrt{1-2kt}\right)$ <br> $= \frac{\left(x - \frac{\mu}{\sqrt{1-2kt}}\right)^2}{2\frac{k}{\sqrt{1-2kt}}x} - \frac{\mu}{k}\left(1 - \sqrt{1-2kt}\right)$ | Dividing the first term by $\sqrt{1-2kt}$ <br> So $t \neq -\frac{1}{2k}$ |

By substitution in the moment expression derive that

$$M_x(t) = e^{\frac{\mu}{k}\left(1-\sqrt{1-2kt}\right)} \int_0^{+\infty} \frac{1}{\sqrt{2\pi k x}} e^{-\frac{\left(x-\frac{\mu}{\sqrt{1-2kt}}\right)^2}{2\frac{k}{\sqrt{1-2kt}}x}} dx$$

$$= \frac{1}{\sqrt{1-2kt}} e^{\frac{\mu}{k}\left(1-\sqrt{1-2kt}\right)} \int_0^{+\infty} \frac{1}{\sqrt{2\pi \frac{k}{1-2kt} x}} e^{-\frac{\left(x-\frac{\mu}{\sqrt{1-2kt}}\right)^2}{2\frac{k}{\sqrt{1-2kt}}x}} dx$$

the last integral is the normalization of an unlikely intensive natural distribution so it equals to 1. So finally

$$M_x(t) = \frac{1}{\sqrt{1-2kt}} e^{\frac{\mu}{k}\left(1-\sqrt{1-2kt}\right)}$$

With analogous proceeding the characteristic function is obtained

$$\Phi_x(t) = \frac{1}{\sqrt{1-2ikt}} e^{\frac{\mu}{k}\left(1-\sqrt{1-2ikt}\right)}.$$

*PROOF 6: Derivation of continuous intensive variable Poisson distribution*

The probability of m events having each a probability p, in n events, when n diverges and p tends to zero as the product Np remain finite, is given by the Poisson distribution

$$pr(m,n) = \frac{1}{m!} np^{m} e^{-np}$$

which is known to have an expectation value E(m) = Np and variance var(m) = Np. Let us now consider that m becomes large and let us use the same approximations of the intensive natural model. The use Stirling approximation for m! leads to

$$pr(m,n) = \frac{1}{\sqrt{2\pi m}} \cdot \left(\frac{np}{m}\right)^{m} \cdot e^{-(np-m)}$$

the extreme value m = 0 can make the probability expression without significance, but as n increase it's probability tends to zero. Then considering that

$$\kappa = \left(\frac{np}{m}\right)^{m} \cdot e^{-(np-m)} \quad \text{and} \quad \log\kappa = -\lambda + m\log\left(1+\frac{\lambda}{m}\right)$$

with λ = np-m. But being m large and λ centered on zero, using the expansion in McLaurin series

$$\begin{aligned}\log\kappa &= -\lambda + m\left\{\left(\frac{\lambda}{m}\right) - \frac{1}{2}\left(\frac{\lambda}{m}\right)^{2} + \frac{1}{3}\left(\frac{\lambda}{m}\right)^{3} - \frac{1}{4}\left(\frac{\lambda}{m}\right)^{4} + \ldots\right\} \\ &= -\frac{1}{2}\frac{\lambda^{2}}{m} + \frac{1}{3}\frac{\lambda^{3}}{m^{2}} - \frac{1}{4}\frac{\lambda^{4}}{m^{3}} + \ldots\end{aligned}$$

In this sum all terms tends to zero as m diverges and the main term remain the first so

$$\log\kappa \approx -\frac{\lambda^{2}}{2m} = -\frac{(Np-m)^{2}}{2m}$$

then the expression is

$$pr(m,n) = \frac{1}{\sqrt{2\pi m}} \cdot e^{-\frac{(Np-m)^{2}}{2m}}$$

Let us then express the probability of the ratio m/n. Consider that

$$\mathrm{pr}\left(\frac{m}{n}\right) = \mathrm{pr}(m)\cdot \mathrm{pr}(n) = \mathrm{pr}(m)$$

In fact we exactly consider n events so pr(n) = 1.

Then $$\mathrm{pr}\left(\frac{m}{n}\right) = \frac{1}{\sqrt{2\pi m}}\cdot e^{-\frac{\left(p-\frac{m}{n}\right)^2}{2\frac{m}{n^2}}}$$

The frequency function remain finite dividing pr(m/n) by the interval 1/n which separates two contiguous values of m/n

$$f\left(\frac{m}{n}\right) = \frac{\mathrm{pr}\left(\frac{m}{n}\right)}{\frac{1}{n}} = \frac{1}{\sqrt{2\pi\frac{m}{n^2}}}\cdot e^{-\frac{\left(p-\frac{m}{n}\right)^2}{2\frac{m}{n^2}}}$$

Now let us measure this very little ratio m/n in an appropriate measurement unity not infinitesimal by multiplication by the factor u. So calling x = um/n, μ = up and k = u/$n^2$ and changing the variable we would obtain

$$f(x) = \frac{1}{\sqrt{2\pi kx}}\cdot e^{-\frac{(\mu-x)^2}{2kx}}$$

which is exactly the intensive natural frequency function for low presence.

So this can be seen as a continuous intensive variable Poisson distributed.